\documentclass[a4paper,10pt,twoside]{amsart}
 \usepackage[latin1]{inputenc}
 \usepackage{graphics}
 \usepackage[dvips]{graphicx}
 \usepackage{graphicx}
\usepackage{subfigure}
 \usepackage{epsfig}
 \usepackage{psfrag}
 \usepackage{amsmath}
 \usepackage{mathrsfs}
 \usepackage{amsfonts}
 \usepackage{amssymb}
 \usepackage{amsthm}
 \DeclareMathOperator{\compo}{o}
 \DeclareMathOperator{\dverg}{div}
 \DeclareMathOperator{\grad}{grad}
 \DeclareMathOperator{\rot}{rot}
 \newcommand{\Tr}{\mathbb{T}}
 \newcommand{\W}{\Omega}
 \newcommand{\w}{\omega}
 
 \newcommand{\eps}{\varepsilon}
 \newcommand{\sob}[2][{}]{{H}^{#1}({#2})}
 \newcommand{\leb}[2][{}]{{L}^{#1}({#2})}
 \newcommand{\V}{\mathbb{V}}
 \newcommand{\U}{\mathbb{U}}
 \newcommand{\N}{\mathbb{N}}
 \newcommand{\Er}{\mathbb{R}}
 \newcommand{\Zr}{\mathbb{Z}}
 
 \newcommand{\Cyl}{\boldsymbol{C}}
 \newcommand{\xdelta}{\mathrm{\delta}}
 \newcommand{\frm}{\mathrm{f}}
 
 \newcommand{\gfrak}{\mathfrak{g}}
\newcommand{\kfrak}{\mathfrak{K}}
 \newcommand{\Cscr}{\mathscr{C}}

 \newcommand{\ri}{\mathrm{i}}
 \newcommand{\xdif}{\mathrm{d}}

 \DeclareMathOperator{\id}{Id}
 \DeclareMathOperator{\inte}{int}
 \DeclareMathOperator{\ext}{ext}
 \DeclareMathOperator{\dvol}{dvol}

 \newtheorem{thrm}{Theorem}
 \newtheorem{nota}[thrm]{Notation}
 \newtheorem{thrm*}{Theorem}
 \newtheorem{lemm}[thrm]{Lemma}
 \newtheorem{rmrk}[thrm]{Remark}
 \newtheorem{rmrk*}{Remark}
 \newtheorem{hypo}[thrm]{Hypothesis}
 
\begin{document}
\title[Asymptotics for Steady State Voltage Potentials]{Asymptotics for Steady State Voltage Potentials in a Bidimensional
  Highly Contrasted Medium with Thin Layer}\thanks{The author thanks Michelle
  Schatzman for many fertile discussions and good advices.}
\author{Clair Poignard}
\address{Centre de Math\'ematiques Appliqu\'ees \\UMR CNRS 7641 \and
  Ecole Polytechnique\\91128 Palaiseau\\ poignard@cmapx.polytechnique.fr}
\date{...}
 \begin{abstract} 
 We study the behavior of steady state voltage potentials in two kinds of  
 bidimensional media composed of material of complex permittivity equal to 1 (respectively $\alpha$)
surrounded by a thin membrane of thickness $h$ and of complex permittivity 
 $\alpha$ (respectively 1).
 We provide in both cases a rigorous derivation of the asymptotic
 expansion of steady state voltage potentials at any order as $h$
 tends to zero, when Neumann boundary condition is imposed on the
 exterior boundary of the thin layer. 
 Our complex parameter $\alpha$ is bounded but may be very small
 compared to 1, hence our results describe the asymptotics of steady
 state 
voltage potentials in all heterogeneous and highly heterogeneous media with thin layer.
 The terms of the potential in the membrane are given explicitly 
 in local coordinates in terms of the boundary data and of the
 curvature of the domain, while these of the inner potential 
are the solutions to the so-called dielectric formulation with 
 appropriate boundary conditions. The error estimates are given explicitly in terms of $h$ and $\alpha$ with appropriate Sobolev 
norm of the boundary data.
We show that the two situations described above lead to completely different asymptotic behaviors of the potentials.

  \end{abstract}
 \maketitle
 \section*{Introduction}
 We study the behavior of the steady state voltage potentials in highly contrasted media surrounded by a thin layer. 
The motivation of the present work comes from numerical problems raised by the researchers in computational electromagnetics, who want to compute the 
 quasi-static electric field in highly contrasted materials with thin layer. 
The thinness of the membrane surrounding an inner domain leads to numerical difficulties, in particular for the meshing. 

To avoid these difficulties, we perform an asymptotic expansion of the potentials in terms of the membrane thickness. 
The approached inner potential is then the finite sum of the solutions to elementary problems in the inner domain with appropriate 
conditions on its boundary, which approximate the effect of the thin layer. 
Thereby, the thin membrane does not have to be considered anymore. Our method leads to the construction of so-called ``approximated 
boundary conditions'' at any order \cite{engquist}. We estimate precisely the error performed by this method in terms of an
appropriate power of 
the relative thinness and with a precise Sobolev norm of the boundary data.
This method is well-known for non highly contrasted media. It is formally described in some particular cases in 
\cite{ammariIEEE} and \cite{idemen}.
We also refer to Kr\"ahenb\"uhl and
 Muller  \cite{krakra} for electromagnetic considerations. Usually, when it is estimated (see for example \cite{engquist}) , the norm of the error involves a imprecise
norm of the boundary data (a $\Cscr^{\infty}$ norm while a weaker norm is enough) and mainly, the constant of the estimate depends 
strongly on the dielectric parameters of the domain.
 
It is 
not obvious (it is even false in general!) that such results hold for highly contrasted domains with thin layer, and this is
 a fact that researchers in computational 
electromagnetics are often confronted to such media. 
For example, a simple electric modelization of the biological cell consists of a conducting cytoplasm surrounded by a thin insulating 
membrane\footnote{We refer to the author thesis~\cite{these} for a precise 
description of the biological cell.}; the modulus of the cytoplasmic complex permittivity divided by the membrane permittivity 
is around $10^5$ while the relative thinness is equal to $10^{-3}$.
On the other hand, the medium might be a dielectric surrounded by a thin metallic layer. 
In both cases it is not clear that the usual approximated boundary 
condition might be used. 

We derive asymptotics 
of the potential steady state voltage in all possible domains with thin layer 
(heterogeneous or highly heterogeneous). As we will see, 
the two situations described above lead to different behaviors of the potentials.
The membrane relative thickness is equal to $h$, while the charateristic
length of the inner domains is equal to 1.
The first medium
consists of a conducting inner domain (say that its complex
permittivity is equal to 1) surrounded by a thin membrane; 
we denote by $\alpha$ the membrane complex permittivity. The parameter
$\alpha$ is bounded but it may tend to zero. This is the reason
why we say that the thin layer is an insulating membrane.
The second material consists of an 
insulating inner domain of permittivity $\alpha$ surrounded by a
conducting thin membrane (say that its complex permittivity is equal
to 1). In this case, we suppose that $\alpha$ tends to zero.
These two kinds of media describe all the possible media with thin
layer.
The aim of this paper is to derive full rigorous asymptotic expansion
of steady state voltage potentials with respect to the small parameter
$h$ for bounded $\alpha$ (but it may tend to zero). 

\begin{figure}[hbt]
 \begin{center}
 \psfrag{a}{$q_h=\alpha\,(\text{resp.}\gamma_h=1)$}
 \psfrag{eps}{$q_h=1(\text{resp.}\gamma_h=\alpha)$}
 \psfrag{ep}{$h$}
 \psfrag{h}{$h$}
 \psfrag{p}{$\mathcal{O}$}
 \psfrag{n}{$\mathcal{O}_{h}$}
 \psfrag{f}{$h$}
 \psfrag{df}{$\W_h$}
 \psfrag{Ce}{${\W}_{l,h}$}
\includegraphics[scale=0.28]{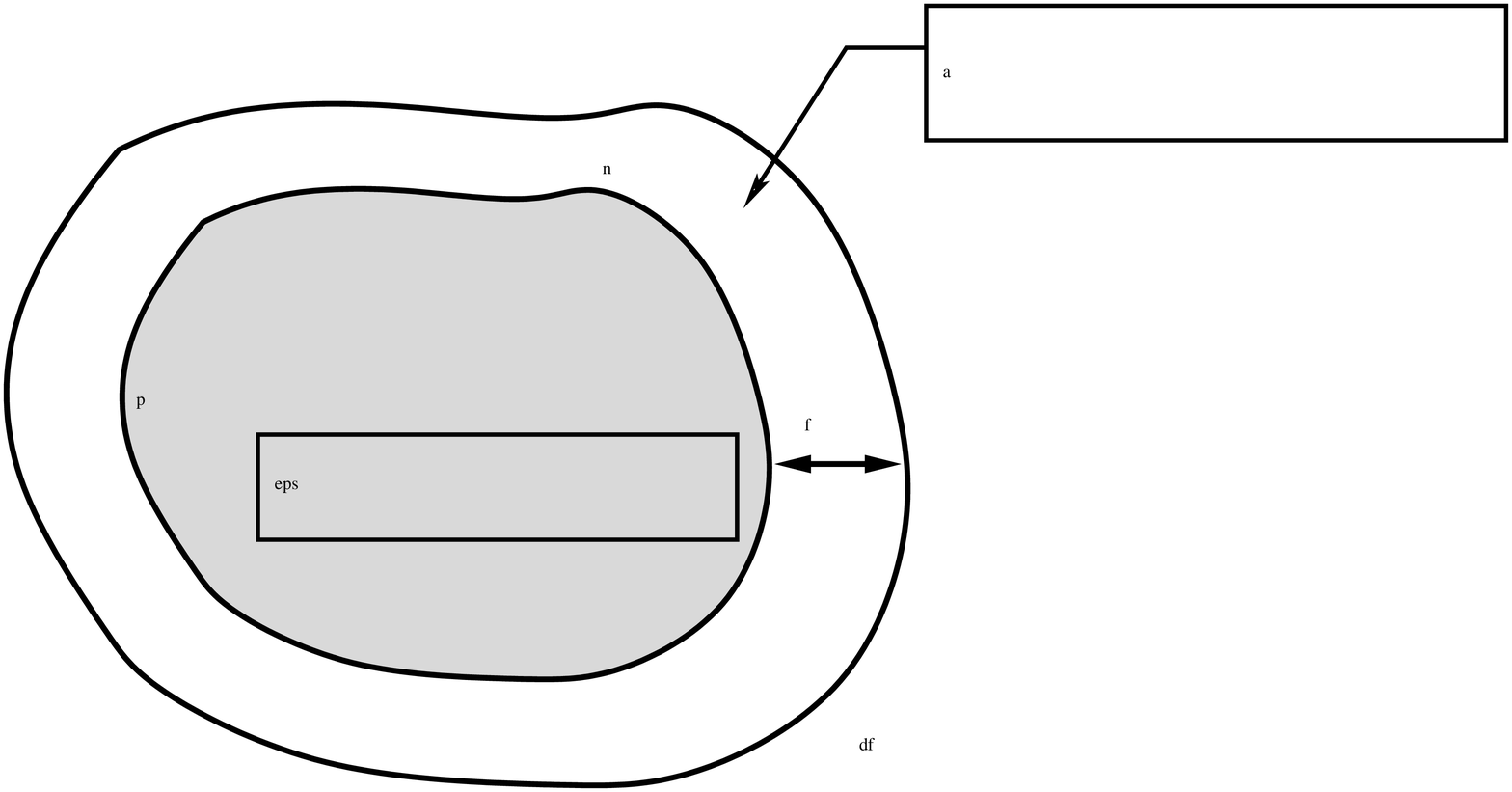}
 \caption{Parameters of $\W_h$.}\label{cell}
  \end{center}
  \end{figure}

Let us write mathematically our problem.
 Let $\W_h$ be a smooth bounded bidimensional 
 domain (see Fig.~\ref{cell}), composed of a smooth domain $\mathcal{O}$ surrounded by a thin membrane $\mathcal{O}_h$ with a 
  small constant thickness~$h$:
 \begin{align*}
 \W_h=\mathcal{O}\cup\mathcal{O}_h.
 \end{align*}
Let $\alpha$ be a non null complex parameter with positive real part; 
$\alpha$ is bounded but it may be very small. Without loss of generality, we suppose that $|\alpha|\leq 1$.
Denote by $q_h$ and $\gamma_h$ the following piecewise constant functions
 \begin{align*}
 \forall x\in\W_h,\quad q_h(x)=\begin{cases}1,\text{ if $x\in \mathcal{O}$},\\
 \alpha,\text{ if $x\in \mathcal{O}_h$},\end{cases}\\
\forall x\in\W_h,\quad \gamma_h(x)=\begin{cases}\alpha,\text{ if $x\in \mathcal{O}$},\\
 1,\text{ if $x\in \mathcal{O}_h$}.\end{cases}
  \end{align*}

  We would like to understand the behavior for $h$  tending to zero
  and uniformly with respect to $|\alpha|\leq1$ of  $V_h$ and $u_h$ 
the respective solutions to the following problems \eqref{diel1} and  \eqref{diel11} with Neumann boundary condition; $V_h$ satisfies
   \begin{subequations}
  \begin{align}
   \nabla\cdot\left(q_h\nabla V_h\right)&=0\label{eqdiel1}
  \text{ in $\W_h$},\\
   \frac{\partial V_h}{\partial n}&=\phi
  \text{ on $\partial\W_h$},\label{cldiel1N}\\
  \int_{\partial \mathcal{O}}V_h\,\xdif \sigma&=0;\label{gauge1}
  \end{align}\label{diel1}
  \end{subequations}
and $u_h$ satisfies
  \begin{subequations}
  \begin{align}
   \nabla\cdot\left(\gamma_h\nabla u_h\right)&=0\label{eqdiel11}
  \text{ in $\W_h$},\\
   \frac{\partial u_h}{\partial n}&=\phi
  \text{ on $\partial\W_h$},\label{cldiel1N1}\\
  \int_{\partial \mathcal{O}}u_h\,\xdif \sigma&=0.\label{gauge11}
  \end{align}\label{diel11}
  \end{subequations}
 Since we impose a Neumann boundary conditions on $\partial\W_h$ the
 boundary data $\phi$ must satisfy the compatibility condition:
 \begin{align*}
 \int_{\partial\W_h}\phi\,\xdif\sigma=0.
 \end{align*}
The above functions $V_h$ and $u_h$ are  well-defined and belong to $\sob[1]{\W_h}$ as soon as $\phi$ belongs to $\sob[-1/2]{\partial\W_h}$. 

Several authors have worked on similar problems (see for instance Beretta \textit{et al.} \cite{beretta2} and  \cite{beretta1}). 
 They compared the exact solution to the so-called background solution
 defined by replacing the material of the membrane by the inner material. 
 The difference between these two solutions has then been given through an 
 integral involving the polarization tensor defined for instance in \cite{ammari1}, \cite{ammari2},  \cite{beretta2},  \cite{beretta1}, \cite{capdebosq}, plus some remainder terms. 
 The remainder terms are estimated in terms of the measure of the inhomogeneity. 
 In this paper, we do not use this approach, for several reasons.

 The Beretta \textit{et al.} estimate of the remainder terms depends
 linearly on $\alpha$ and  $1/\alpha$: their results are no more valid
 in a highly contrasted domain (\textit{i.e.} for $\alpha$ very large
 or very small).
 Secondly, $\alpha$ is complex-valued, hence differential operators involved in our case are not self-adjoint, so that the $\Gamma$-convergence techniques of Beretta \textit{et al.} do not apply.
 Thirdly, the potential in the membrane is not given explicitly in \cite{beretta2}, \cite{beretta1} or \cite{capdebosq}, while we are definitely interested in this potential, in order to obtain 
 the transmembranar potential (see Fear and Stuchly \cite{fear}). Finally, the asymptotics of Beretta \textit{et al.} are valid on the boundary of the domain, while we are interested in the potentials in the inner domain.

 The heuristics of this work consist in performing a change of coordinates in the membrane $\mathcal{O}_h$, so as to parameterize it by local coordinates $(\eta,\theta)$, which vary in a domain 
 independently 
 of $h$; in particular, if we denote by $L$ the length of $\partial
 \mathcal{O}$ (in the following, without any restriction, we suppose that $L$ is equal to $2\pi$)), the variables $(\eta,\theta)$ should vary in
 $[0,1]\times \Er/L\Zr$. This change of coordinates leads to an expression of
 the Laplacian in 
 the membrane, which depends on $h$. Once the transmission conditions of the new problem are derived, we perform a formal asymptotic expansion of the solution to Problem~\eqref{diel1} (respectively to Problem \eqref{diel11}) in terms of $h$. It remains to
  validate this expansion. In this paper we work with bidimensional domain and we are confident 
that the same analysis could be perfomed in higher dimensions. 

 This paper is structured as follows. In Section~\ref{geo}, we make
 precise our geometric conventions. We perform a change of variables
 in the membrane, 
and with the help of some differential geometry results, 
 we write Problem \eqref{diel1}  and Problem~\eqref{diel11} in the language of differential forms. 
 We refer the reader to Flanders \cite{flanders} or Dubrovin
 \textit{et al.} \cite{dubrovin} (or \cite{doubrovine} for the french
 version) for courses on 
differential geometry. 
 We derive transmission and boundary conditions in the intrinsic
 language of differential forms, and we express these relations in
 local coordinates.

In Section~\ref{insulating membrane} we study Problem~\eqref{diel1}. In paragraph~\ref{formal}  we derive formally all the terms
 of the  asymptotic expansion of the solution to our problem in terms
 of $h$. Paragraph~\ref{estimates} is devoted to a proof of the
 estimate of the error. 

Problem~\eqref{diel11} is
considered in Section~\ref{insulating inner domain}. We supposed that $\alpha$ tends to zero: a boundary layer
phenomenon appears. To obtain our error estimates, we link the
parameters $h$
and $\alpha$. We introduce a complex parameter $\beta$ such that
$$Re(\beta)>0, or \left(Re(\beta)=0,\text{ and } \Im(\beta)\neq0\right),$$
and
$$|\beta|=o\left(\frac{1}{h}\right),\quad\text{and}\quad \frac{1}{|\beta|}=o\left(\frac{1}{h}\right).$$
We distinguish two different cases, depending on the convergence 
of $|\alpha|$ to zero: $\alpha=\beta h^q$,  for $q\in N^*$ and $\alpha=o(h^{N})$ for all $N\in\N$.

For $q=1$ we obtain mixed boundary conditions for the asymptotic terms
of the inner potential, and as soon as $q\geq1$, appropriate Dirichlet
boundary conditions are obtained. We end this section by error estimates.

In Appendix, we give some useful differential
  geometry  formulae.

\begin{rmrk*}
The use of the formalism of differential forms  $\delta\left(q_h\xdif \right)$ could seem futile for the study of the operator $\nabla\cdot\left(q_h \nabla\right)$. In particular the expression of Laplace 
operator in local coordinates is well known. 
However we wanted to present this point of view to show how simple it is to write a Laplacian in curved coordinates once the metric is known.

Moreover once this formalism is understood for the functions (or 0-forms), it is easy to study $\delta\left(q_h\xdif \right)$ applied to 1-forms. This leads directly to the study of the operator
 $\rot\left(q_h\rot \right),$ whose expression in local coordinates is less usual.
\end{rmrk*}
We choose to present our two main theorems in this introduction so that the reader interested in our results 
without their proves might find them easily.

We suppose that
 $\partial{\mathcal{O}}$ is smooth. 
We denote by $\Phi$ the
$\Cscr^\infty-$ diffeormorphism,
which maps a neighborhood of cylinder $\Cyl=[0,1]\times\Er/2\pi\Zr$ unto a
neighborhood of the thin layer. The diffeomorphism
$\Phi_0=\Phi(0,\cdot)$ maps the torus unto the boundary $\partial
\mathcal{O}$ of the inner domain while $\Phi_1=\Phi(1\cdot)$ is the
$\Cscr^\infty-$diffeomorphism from the torus unto $\partial \W_h$. We
denote by $\kappa$ the curvature of $\partial\mathcal{O}$ written in
local coordinates, and let $h_0$ be such that
$$h_0<\frac{1}{\sup_{\theta\in\Er/2\pi\Zr}|\kappa(\theta)|}.$$

\subsection*{Asymptotic for an insulating thin layer}
The first theorem gives the
asymptotic expansion of the solution $V_h$ of \eqref{diel1}, for $h$ tending to zero,
for bounded $\alpha$.
\begin{thrm}\label{thm1}
Let $h$ belong to $(0,h_0)$.
 The complex parameter $\alpha$  satisfies
 \begin{align}
&|\alpha|\leq 1,\label{hyp1}\\
&\Re(\alpha)> 0\,\text{ or }\,\Bigl\{\Re(\alpha)=0\,\text{and }\,\Im(\alpha)\neq0\Bigr\}.\label{hyp}
 \end{align}
 Let $N\in\N$ and $\phi$ belong to $\sob[N+3/2]{\partial \W_h}$.  Denote by $f$ and $\frm$ the following functions:
\begin{align*}
\forall \theta\in\Er/2\pi\Zr,\quad f(\theta)&=\phi\compo\Phi_1(\theta),\\
\forall x\in\partial\mathcal{O},\quad \frm(x)&=\phi\compo\Phi_1\compo\Phi^{-1}_0(x).\end{align*}
Define the sequence of potentials $(V^c_k,V^m_k)^N_{k=0}$ as follows.
We impose
\begin{align*}
&\forall (\eta,\theta)\in\Cyl,\quad \partial_\eta V^m_0=0,
\end{align*}
and we use the convention 
\begin{align*}\begin{cases}
V^{c}_l=0,\text{ if $l\leq -1$},\\
V^{m}_l=0,\text{ if $l\leq -1$}.\end{cases}
\end{align*} 
For $0\leq k\leq N$ we define for all $0\leq s\leq 1$ the function
$\partial_\eta V^m_{k}(s,\cdot)$ on $\Er/2\pi\Zr$ :
\begin{align*}
\begin{split}\partial_\eta V^m_{k+1}(s,\cdot)=&\delta_{1,k+1}f+\int^1_s\biggl\{\kappa\left\{3\eta\partial^2_\eta V^m_{k}+\partial_\eta V^m_{k}\right\}
\\&+3\eta^2\kappa^2\partial^2_\eta V^m_{k-1}+2\eta\kappa^2\partial_\eta V^m_{k-1}+\partial^2_\theta V^m_{k-1}
\\&+\eta^3\kappa^3\partial^2_\eta V^m_{k-2}+\eta^2\kappa^3\partial_\eta V^m_{k-2}+\eta\kappa\partial^2_\theta V^m_{k-2}
-\eta\kappa'\partial_\theta V^m_{k-2}\biggr\}\xdif \eta,\end{split}
\end{align*}
and  the functions $V^c_k$ and $V^m_k$ are then defined by
\begin{align*}
&\left.\Delta V^c_{k}\right.=0,\\
&\left.\partial_n V^c_{k}\right|_{\partial\mathcal{O}}=\alpha \partial_\eta V^m_{k+1}\compo\Phi^{-1}_0,\\
&\int_{\partial\mathcal{O}}V^c_{k}\xdif \sigma=0,\\
&\forall s\in(0,1),\quad V^m_{k}(s,\cdot)=\int^s_0 \partial_\eta
  V^m_{k}(\eta,\cdot)\,\xdif\eta+V^c_{k}\compo\Phi_0.\
\end{align*}
Let $R^c_N$ and $R^m_N$ be the functions defined by:
 \begin{align*}
\begin{cases} 
R^c_N=V_h-\sum^N_{k=0}V^c_kh^k,\,\text{in $\mathcal{O}$},\\
R^m_N=V_h\compo\Phi-\sum^N_{k=0}V^m_kh^k,\,\text{in $\Cyl$}.
 \end{cases}
 \end{align*}
 Then, there exists a constant $C_{\mathcal{O},N}>0$ depending only on the domain $\mathcal{O}$ and on $N$ such that 
\begin{subequations}
 \begin{align}
\|R^c_N\|_{\sob[1]{\mathcal{O}}}&\leq C_{\mathcal{O},N}\|\frm\|_{\sob[N+3/2]{\partial\mathcal{O}}}|\alpha| h^{N+1/2},\\
\|R^m_N\|_{H^1_\gfrak\left(\Cyl\right)}&\leq C_{\mathcal{O},N}\|\frm\|_{\sob[N+3/2]{\partial\mathcal{O}}}h^{N+1/2}.
 \end{align}\label{resultat11}
\end{subequations}
Moreover, if $\phi$ belongs to $\sob[N+5/2]{\partial \W_h}$, then we have
\begin{subequations}
 \begin{align}
\|R^c_N\|_{\sob[1]{\mathcal{O}}}&\leq C_{\mathcal{O},N}\|\frm\|_{\sob[N+5/2]{\partial\mathcal{O}}}|\alpha| h^{N+1},\\
\|R^m_N\|_{H^1_\gfrak\left(\Cyl\right)}&\leq C_{\mathcal{O},N}\|\frm\|_{\sob[N+5/2]{\partial\mathcal{O}}}h^{N+1/2}.
 \end{align}\label{resultat21}
\end{subequations}
 \end{thrm}
In this theorem, we approach the potential in the inner domain at the
order $N$ by solving $N$ elementary problems with appropriate boundary
condition. From these results, we may build an approximated boundary condition
on $\partial\mathcal{O}$
at any order, in order to solve only one problem. However, this kind of
conditions lead to numerical unstabilities, this is the reason why we
think that the method to obtain the potential step by step is more
useful. 

Since it is classical to write approximated boundary conditions
we make precise these conditions at the orders 0 and 1. 
Denote by $\kfrak$ the curvature of $\partial\mathcal{O}$ in Euclidean coordinates
and by
$V^0_{app}$  and $V^1_{app}$ the approximated potentials with
approximated boundary condition at the order 0 and 1 respectively. We
have:
\begin{align}
&\Delta V^0_{app}=0,\text{ in $\mathcal{O}$},\\
&\partial_nV^0_{app}=\alpha \frm,\label{CLapp0}
\intertext{and }
&\Delta V^1_{app}=0,\text{ in $\mathcal{O}$},\\
&\partial_nV^1_{app}-\alpha h \partial^2_t V^1_{app}=\alpha (1+h\kfrak)\frm,\label{CLapp1}
\end{align}
where $\partial_t$ denotes the tangential derivative on $\partial\mathcal{O}$.
The boundary condition \eqref{CLapp1} imposed to $\partial_nV^1_{app}$ is well-known
for non highly contrasted media. It might be found in \cite{krakra}. With our theorem, we prove that it remains valid for a
very insulating membrane, and we give precise norm
estimates. Moreover we give complete asymptotic expansion of the
potential in both domains (the inner domain and the thin layer). 

\begin{figure}[h!]
\begin{center}
\subfigure[$\alpha= i$]{
\includegraphics[angle=0, scale=0.6]{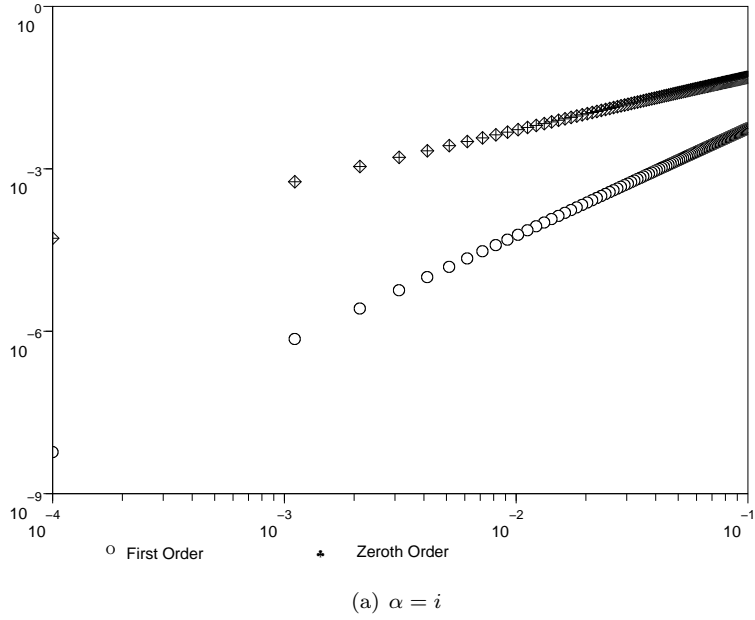}}\label{a}
\subfigure[$\alpha=i/h$]{
\includegraphics[angle=0, scale=0.6]{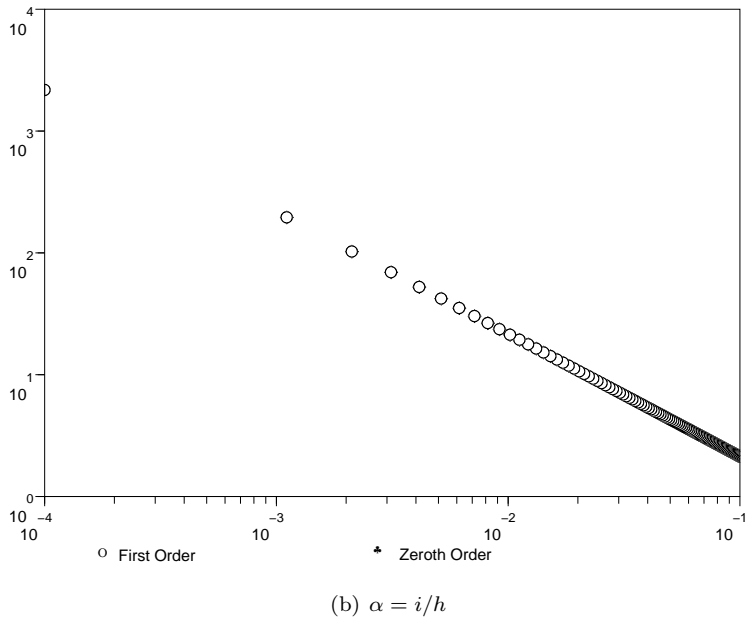}}\label{b}
 \caption{$H^1$ norm of the error at the orders 0 and 1.}\label{CL10usuel}
\end{center}  
\end{figure}

We perform numerial simulations in a circle of radius 1 surrounded by a thin layer of thickness $h$.
In Fig~\ref{CL10usuel}, the left frame illustrates the asymptotic estimates at the orders 0 and 1 of Theorem~\ref{thm1} 
for an insulating thin layer. However, the right frame shows that as soon as the thin layer becomes very conducting, 
for example as soon as $\alpha=1/h$, 
these asymptotics are no more valid: we have to use the asymptotics of Theorem~\ref{thm2}.

\subsection*{Asymptotics for an insulating inner domain}
\begin{figure}[hbt]
 \begin{center}
\includegraphics[angle=0, scale=0.6]{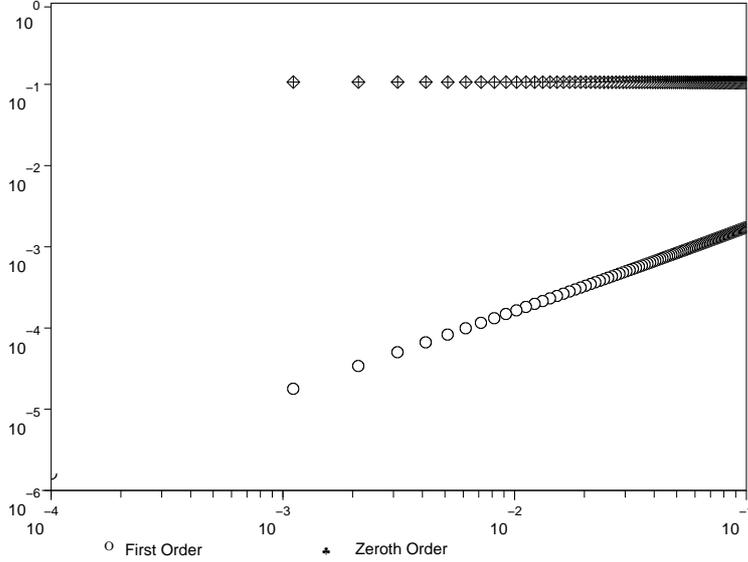}
 \caption{$H^1$ Norm of the error at the order $-1$ and $0$ for an insulating inner domain : $\alpha=i h$.}\label{fig2}
  \end{center}
  \end{figure}

Let $\beta$ be a complex parameter satisfying:
$$Re(\beta)>0, or \left(Re(\beta)=0,\text{ and } \Im(\beta)\neq0\right).$$
The modulus of $\beta$ may tend to infinity, or to zero but it must satisfy:
$$|\beta|=o\left(\frac{1}{h}\right),\quad\text{and}\quad \frac{1}{|\beta|}=o\left(\frac{1}{h}\right).$$
 \begin{thrm}\label{thm2}
 Let $h$ belong to $(0,h_0)$. Let $q\in\N^*$ and $N\in\N$.
We suppose that  $\alpha$ satisfies:
\begin{align}
\alpha&=\beta h^q.
\end{align} 
Let $\phi$ belong to $\sob[N+3/2+q]{\partial \W_h}$ and denote by $f$ and $\frm$ the following functions:
\begin{align*}
\forall \theta\in\Er/2\pi\Zr,\quad f(\theta)&=\phi\compo\Phi_1(\theta),\\
\forall x\in\partial\mathcal{O},\quad \frm(x)&=\phi\compo\Phi_1\compo\Phi^{-1}_0(x).
\end{align*}
Define the function $(u^{c,q}_k,u^{m,q}_k)^N_{k=-1}$ by induction as
 follows, with the convention 
\begin{align*}\begin{cases}
u^{c,q}_l=0,\text{ if $l\leq -2$},\\
u^{m,q}_l=0,\text{ if $l\leq -2$}.\end{cases}
\end{align*} 
\begin{itemize}
\item If $q=1$
\begin{equation*}
\left\{\begin{aligned}
&\Delta u^{c,1}_{-1}=0,\text{ in $\mathcal{O}$,}\\
&-\left.\partial^2_t u^{c,1}_{-1}\right|_{\partial\mathcal{O}} +\beta\left.\partial_n u^{c,1}_{-1}\right|_{\partial\mathcal{O}}=\frm,\\
&\int_{\partial\mathcal{O}}u^{c,1}_{-1}\xdif_{\partial\mathcal{O}}=0.
\end{aligned}\right.
\end{equation*}
\begin{equation*}
\forall (\eta,\theta)\in\Cyl,\quad u^{m,1}_{-1}=u^{c,1}_{-1}|_{\partial\mathcal{O}}\compo\Phi_0.
\end{equation*}
Moreover, \begin{align*}&\partial_\eta u^{m,1}_{0}=0,\quad 
\partial_\eta u^{m,1}_{1}=(1-\eta)\partial^2_\theta u^{m,1}_{-1}+f.\end{align*} 
For  $0\leq k\leq N$, denote by $\phi^1_{k}$ the following function:
$$\phi^1_k=\int^1_0\Bigl(\kappa\left(3\eta\partial^2_\eta u^{m,1}_{k+1}+\partial_\eta u^{m,1}_{k+1}\right)+\eta\kappa\partial^2_\theta u^{m,1}_{k-1}-\eta\kappa'
\partial_\theta u^{m,1}_{k-1}\Bigr)\,\xdif\eta.$$
and define $u^{c,1}_{k}$ by
\begin{equation*}
\left\{\begin{aligned}
&\Delta u^{c,1}_{k}=0,\text{ in $\mathcal{O}$,}\\
&-\left.\partial^2_t u^{c,1}_{k}\right|_{\partial\mathcal{O}} +\beta\left.\partial_n u^{c,1}_{k}\right|_{\partial\mathcal{O}}
=\Biggl(\phi^1_k-\int^1_0(\eta-1)\partial^2_\theta\partial_\eta u^{m,1}_k\xdif \eta\Biggr)\compo\Phi^{-1}_0,\\
&\int_{\partial\mathcal{O}}u^{c,1}_{k}\xdif_{\partial\mathcal{O}}=0.
\end{aligned}\right.
\end{equation*}
{In the membrane $u^{m,1}_k$ is defined by }
\begin{equation*}
u^{m,1}_k=\int^s_0\partial_\eta u^{m,q}_{k}\xdif\eta+u^{c,q}_k\compo\Phi_0,
\end{equation*}
and $\partial_\eta u^{m,1}_{k+i}$ for ${i=1,2}$ is determined by:
\begin{align}
\begin{split}\partial_\eta u^{m,1}_{k+i}&=\int^s_1\Biggl(-\kappa\left(3\eta\partial^2_\eta u^{m,1}_{k+i-1}
+\partial_\eta u^{m,1}_{k+i-1}\right)\\&-\partial^2_\theta u^{m,1}_{k+i-2}-\eta\kappa\partial^2_\theta u^{m,1}_{k+i-3}+\eta\kappa'
\partial_\theta u^{m,1}_{k+i-3}\Biggr)\xdif\eta.\end{split}\label{dnumkp1}\end{align}
\item If $q\geq 2$.
The function $u^{m,q}_1$ is defined by
\begin{align*}
&\int_\Tr u^{m,q}_{-1}\xdif\theta=0,\\
&-\partial^2_\theta u^{m,q}_{-1}=f.
\end{align*}
The potential $u^{c,q}_{-1}$ is solution to the following problem:
\begin{equation*}
\left\{\begin{aligned}
&\Delta u^{c,q}_{-1}=0,\text{ in $\mathcal{O}$,}\\
&\left. u^{c,q}_{-1}\right|_{\partial\mathcal{O}}=u^{m,q}_{-1}\compo\Phi^{-1}_0 .
\end{aligned}\right.
\end{equation*}
Moreover, \begin{align*}&\partial_\eta u^{m,q}_{0}=0,\quad 
\partial_\eta u^{m,q}_{1}=(1-\eta)\partial^2_\theta u^{m,q}_{-1}+f.\end{align*} 
For $0\leq k\leq N$, denote by $\phi^q_{k}$ the following function:
$$\phi^q_k=\int^1_0\Bigl(\kappa\left(3\eta\partial^2_\eta u^{m,q}_{k+1}+\partial_\eta u^{m,q}_{k+1}\right)+\eta\kappa\partial^2_\theta u^{m,q}_{k-1}-\eta\kappa'
\partial_\theta u^{m,q}_{k-1}\Bigr)\,\xdif\eta.$$
$u^{m,q}_k|_{\eta=1}$ is entirely determined by the equality:
\begin{align*}-\partial^2_\theta u^{m,q}_k|_{\eta=1} 
=&\beta\partial_n u^{c,q}_{k+1-q}\compo\Phi_0+\phi^q_k-\int^1_0\eta\partial^2_\theta\partial_\eta u^{m,q}_k\xdif \eta,\intertext{hence}
u^{m,q}_k(s,\theta)&=\int^s_1\partial_\eta u^{m,q}_{k}\xdif\eta+u^{m,q}_k|_{\eta=1}.
\end{align*}
The potential $u^{c,q}_k$ satisfies the following boundary value problem:
\begin{equation*}
\left\{\begin{aligned}
&\Delta u^{c,q}_{k}=0,\text{ in $\mathcal{O}$,}\\
&\left. u^{c,q}_{k}\right|_{\partial\mathcal{O}}=u^{m,q}_{k}\compo\Phi^{-1}_0 .
\end{aligned}\right.
\end{equation*}
The functions $\left(\partial_\eta u^{m,q}_{k+i}\right)_{i=1,2}$ satisfies equation \eqref{dnumkp1}, in which $u^{m,1}$ is replaced by $u^{m,q}$.
\end{itemize}
Let $r^{c,q}_N$ and $r^{m,q}_N$ be the functions defined by:
 \begin{align*}
\begin{cases} r^{c,q}_N=u_h-\sum^N_{k=-1}u^{c,q}_kh^k,\,\text{in $\mathcal{O}$},\\
 r^{m,q}_N=u_h\compo\Phi-\sum^N_{k=-1}u^{m,q}_kh^k,\,\text{in $\Cyl$}.
 \end{cases}
 \end{align*}
 Then, there exists a constant $C_{\mathcal{O},N}>0$ depending only on the domain $\mathcal{O}$ and on $N$ such that 
 \begin{align*}
\|r^{c,q}_N\|_{\sob[1]{\mathcal{O}}}&\leq C_{\mathcal{O},N}\|\frm\|_{\sob[N+3/2+q]{\partial\mathcal{O}}} \max\left(\sqrt{\frac{h}{|\beta|}},\sqrt{h}\right)h^{N+1/2},\\
\|r^{m,q}_N\|_{H^1_\gfrak\left(\Cyl\right)}&\leq C_{\mathcal{O},N}\|\frm\|_{\sob[N+3/2]{\partial\mathcal{O}}}h^{N+1/2}.
 \end{align*}
If $\phi$ belongs to $\sob[N+5/2+q]{\partial\W_h}$, we have
$$\left\| r^{c,q}_N\right\|_{\sob[1]{\mathcal{O}}}\leq
 C_{\mathcal{O},N}\|f\|_{\sob[N+5/2+q]{\Tr}}h^{N+1}.$$
 \end{thrm}
We observe that if $q=1$ and $N=0$, the approximated boundary condition
at the order 0 is given by:
$$-(1-h\kfrak/2)\partial^2_tu^{1}_{0,app}+\frac{h\partial_t\kfrak}{2}\partial_tu^{1}_{0,app}+\beta 
\partial_nu^{1}_{0,app}=\frac{1+h\kfrak}{h}\phi\compo\Phi_1\compo\Phi^{-1}_0.$$
Thus it is very different from the approximated boundary condition \eqref{CLapp0} imposed to $V^0_{app}$ in the case of 
an insulating membrane. This is a feature of the conducting thin layer. 
Observe on Fig~\ref{fig2} that the numerical computations 
in a circle confirm our theorical results.

Thanks to our previous results  by comparing the parameters $|\alpha|$ and $h$ of a heterogeneous medium with thin layer, 
we know \textit{a priori}, which asymptotic formula 
(Theorem~\ref{thm1} or Theorem~\ref{thm2}) has to be computed. 
We emphasize that our method might be easily implemented by iterative process as soon as the geometry of the domain 
is precisely known.

In the following, we show how the potentials $(V^c_k, V^m_k)_{k\geq 0}$ of the previous theorems 
are built, and then we prove these theorems. Let us now make precise the geometric conventions.
\section{Geometry}\label{geo}
The boundary of the domain $\mathcal{O}$ is assumed to be smooth. 
The orientation of the boundary $\partial \mathcal{O}$ is the trigonometric orientation. 
To simplify, we suppose that the length of
$\partial \mathcal{O}$ is equal to $2\pi$. We denote by $\Tr$ the flat torus:
\begin{align*}
\Tr=\Er/2\pi\Zr.
\end{align*} 
Since $\partial \mathcal{O}$ is smooth, we can parameterize it
by a function $\Psi$ of class $\Cscr^\infty$ from $\Tr$ to $\Er^2$ satisfying:  
\begin{align*}
\forall \theta\in\Tr,\quad\left|\Psi'\left(\theta\right)\right|=1.
\end{align*}
Since the boundary $\partial \W_h$ of the
cell is parallel to the boundary $\partial \mathcal{O}$ of the inner domain
the following identities hold: 
\begin{align*}
\partial \mathcal{O}&=\{\Psi(\theta),\theta\in\Tr\},
\intertext{and}
\partial \W_h&=\{\Psi(\theta)+hn(\theta),\theta\in\Tr\}.
\intertext{Here $n(\theta)$ is the unitary exterior normal at $\Psi(\theta)$ to $\partial\mathcal{O}$. 
Therefore the membrane $\mathcal{O}_h$ is parameterized by:}
\mathcal{O}_h&=\{\Phi(\eta,\theta),\,
(\eta,\,\theta)\in]0,\,1[\times\Tr\},
\intertext{where}
\Phi(\eta,\theta)&=\Psi(\theta)+h\eta n(\theta).
\end{align*}

Denote by $\kappa$ the curvature of $\partial \mathcal{O}$.
Let $h_0$ belong to $(0,1)$ such that:
\begin{align}
h_0<\frac{1}{\|\kappa\|_\infty}.\label{h0}
\end{align} 
Thus for all $h$ in $[0,h_0]$, there exists an open intervall $I$ containing $(0,1)$ such that $\Phi $ is a smooth diffeomorphism 
from $I\times \Er/2\pi\Zr$ to its image, which is a neighborhood of
the membrane.
The metric in $\mathcal{O}_h$ is:
\begin{align}
h^2\xdif \eta^2+(1+h\eta\kappa)^2\xdif \theta^2.\label{localmetr}
\end{align}
Thus, we use two systems of coordinates, depending on the domains $\mathcal{O}$ and $\mathcal{O}_h$: in the interior domain $\mathcal{O}$, we use Euclidean coordinates
$(x,y)$ and in the membrane $\mathcal{O}_h$, we use local $(\eta,\theta)$ coordinates  with metric \eqref{localmetr}.

We translate  into the language of differential forms Problem
\eqref{diel1} and Problem \eqref{diel11}. We refer the reader to
Dubrovin, Fomenko and Novikov \cite{dubrovin} or 
Flanders \cite{flanders} for the definition of the exterior derivative
denoted by $\xdif$, the exterior product denoted by $\ext$, the
interior derivative denoted by $\delta$ and the interior product
denoted by $\inte$. In Appendix 
we give the formulae describing these operators in 
the case of a 
general 2D metric.
Our aim, while rewriting our problems \eqref{diel1} and \eqref{diel11} is to take into account nicely the change of coordinates in the thin membrane.

Let $\V$ be the 0-form on $\W_h$ such that, in the Euclidean coordinates $(x,y)$, $\V$ is equal to $V$, and let $F$ be the 0-form, which is equal  to $\phi$ on $\partial \W_h$. We denote by $\N$ the 1-form corresponding to the 
inward unit normal on the boundary $\W_h$ (see for instance Gilkey \textit{et al.} \cite{gilkey2} p.33):
\begin{align*}
\N&=N_x\xdif x+N_y \xdif y,\\
&=N_\eta\xdif \eta. 
\end{align*}
$\N^*$ is the inward unit normal 1-form.
Problem \eqref{diel1} takes now the intrinsic form:
\begin{subequations}
\begin{align}
\delta\left(q_h\xdif\V\right)&=0,\,\text{in $\W_h$},\\
\inte(\N^*)\xdif \V&=F,\,\text{on $\partial \W_h$}.\label{cldiff}
\intertext{According to Green's formula (Lemma 1.5.1 of \cite{gilkey2}), we obtain the following transmission conditions for $\V$ along $\partial\mathcal{O}$:}
\begin{split}\inte(\N^*)\xdif \V|_{\partial \mathcal{O}}&=\alpha\inte(\N^*)\xdif \V|_{\partial \mathcal{O}_h\setminus\partial\W_h},\\
\ext(\N^*) \V|_{\partial \mathcal{O}}&=\ext(\N^*) \V|_{\partial \mathcal{O}_h\setminus\partial\W_h}.
\end{split}\label{transmi12}
\end{align}\label{dielfdif}
\end{subequations}
Similarly, denoting by $\U$ the 0-form equal to $u$ in Euclidean
coordinates  we rewrite Problem \eqref{diel11} as follows:
\begin{subequations}
\begin{align}
\delta\left(\gamma_h\xdif\U\right)&=0,\,\text{in $\W_h$},\\
\inte(\N^*)\xdif \U&=F,\,\text{on $\partial \W_h$};\label{cldiffU}
\intertext{the following transmission conditions hold on $\partial\mathcal{O}$:}
\begin{split}\alpha\inte(\N^*)\xdif \U|_{\partial \mathcal{O}}&=\inte(\N^*)\xdif \U|_{\partial \mathcal{O}_h\setminus\partial\W_h},\\
\ext(\N^*) \U|_{\partial \mathcal{O}}&=\ext(\N^*) \U|_{\partial \mathcal{O}_h\setminus\partial\W_h}.
\end{split}\label{transmi12U}
\end{align}\label{dielfdifU}
\end{subequations}

\section{Statement of the problem}\label{pb}
In this section, we write Problem~\eqref{dielfdif} and
Problem~\eqref{dielfdifU} in local coordinates, with the help of
differential forms.
It is convenient to write:
\begin{align*}
\forall \theta\in\Tr,\quad \Phi_0\left(\theta\right)&=\Phi\left(0,\theta\right),\,\Phi_1\left(\theta\right)=\Phi\left(1,\theta\right),
\intertext{and to denote by $\Cyl$ the cylinder:}
\Cyl&=[0,1]\times\Tr.
\end{align*}
We denote by $\kfrak$, $f$ and $\frm$ the following functions:
\begin{align}
\forall (x,y)\in\partial\mathcal{O},\quad \kfrak(x,y)&=\kappa\compo\Phi^{-1}_0(x,y),\\
\forall \theta\in\Tr,\quad f(\theta)&=\phi\compo\Phi_1(\theta),\label{ftangentielle}\\
\forall x\in\partial\mathcal{O},\quad\frm&=f\compo\Phi_0^{-1}(x).\label{frm}
\end{align}
Using the 
expressions of the differential operators $\xdif$ and $\xdelta$, which are respectively the exterior and the interior derivatives (see Appendix), applied to the metric 
\eqref{localmetr}, the Laplacian in the membrane is given in the local coordinates $(\eta,\theta)$ by:
\begin{align}
&\forall\,(\eta,\theta)\in\Cyl,\notag\\ 
\left.\Delta \right|_{\Phi\left(\eta,\theta\right)}=&\frac{1}{h(1+h\eta\kappa)}\partial_\eta\left(\frac{1+h\eta\kappa}{h}\partial_\eta \right)
+\frac{1}{1+h\eta\kappa}\partial_\theta\left(\frac{1}{1+h\eta\kappa}\partial_\theta \right).\label{laplloc}
\end{align}
Moreover, for a 0-form $z$ defined in $\mathcal{O}_h$, we have:
\begin{align*}
\left.\inte(\N^*)\xdif z\right|_{\partial\mathcal{O}}&=\frac{1}{h}\partial_\eta z|_{\eta=0},\\ 
\left.\inte(\N^*)\xdif z\right|_{\partial\W_h}&
=\frac{1}{h}\partial_\eta z|_{\eta=1}.
\end{align*}
Denote by
\begin{align*}
 V^c&=V,\,\text{in $\mathcal{O}$},\\
 V^m&=V\compo\Phi,\,\text{in $\Cyl$},
\intertext{and by}
u^c&=u,\,\text{in $\mathcal{O}$},\\
u^m&=u\compo\Phi,\,\text{in $\Cyl$}.
\end{align*}
We infer that Problem \eqref{dielfdif} may be rewritten as follows:
\begin{subequations}
\begin{align}
&\Delta V^c=0,\,\text{in $\mathcal{O}$},\label{laplVc}\\
\forall\,(\eta,\theta)\in\Cyl,\quad
&\frac{1}{h^2}\partial_\eta\left((1+h\eta\kappa)\partial_\eta
V^m\right)
+\partial_\theta\left(\frac{1}{1+h\eta\kappa}\partial_\theta V^m\right)=0,\label{laplacemetr}\\
&{\partial_n V^c}\compo \Phi_0=\left.\frac{\alpha}{h}{\partial_\eta V^m}\right|_{\eta=0},\label{transmi1}\\
&V^c\compo\Phi_0=\left.V^m\right|_{\eta=0},\label{transmi2}\\
&\left.\partial_\eta V^m\right|_{\eta=1}= {h}f.\label{bcdielmetr}\\
& \int_{\partial \mathcal{O}}V\,\xdif \sigma=0.\notag
\end{align}\label{dielmetr}
\end{subequations}
Similarly the couple $(u^c,u^m)$ satisfies
\begin{subequations}
\begin{align}
&\Delta u^c=0,\,\text{in $\mathcal{O}$},\label{lapluc}\\
\forall\,(\eta,\theta)\in\Cyl,\quad &
\frac{1}{h^2}\partial_\eta\biggl((1+h\eta\kappa)\partial_\eta u^m\biggr)+\partial_\theta\left(\frac{1}{1+h\eta\kappa}\partial_\theta u^m\right)=0,\label{laplacemetru}\\
&\alpha{\partial_n u^c}\compo \Phi_0=\left.\frac{1}{h}{\partial_\eta u^m}\right|_{\eta=0},\label{transmi1u}\\
&u^c\compo\Phi_0=\left.u^m\right|_{\eta=0},\label{transmi2u}\\
&\left.\partial_\eta u^m\right|_{\eta=1}= {h}f,\label{bcdielmetru}\\
& \int_{\partial \mathcal{O}}u\,\xdif \sigma=0.\notag
\end{align}\label{dielmetru}
\end{subequations}
 \begin{rmrk}
In the following, the parameter $\alpha$ is such that:
 \begin{align*}
&\Re(\alpha)> 0\,\text{ or }\,\Bigl\{\Re(\alpha)=0\,\text{and }\,\Im(\alpha)\neq0\Bigr\}.
 \end{align*}
 Since $\alpha$ represents a complex permittivity it may be written
 (see Balanis and Constantine \cite{balanis}) as follows:
 \begin{align*}
\alpha=\eps-\ri \sigma/\w,
 \end{align*}
 with $\eps$, $\sigma$,  and $\w$ positive. Thus this hypothesis is always satisfied for dielectric materials.
 \end{rmrk}

\begin{nota}\label{notasob}
We provide $\Cyl$ with the metric \eqref{localmetr}. The $L^2$ norm of a $0$-form $u$ in $\Cyl$, denoted by $\|u\|_{\Lambda^0L^2_m(\Cyl)}$, is equal to:
 \begin{align*}
 \|u\|_{\Lambda^0L^2_\gfrak\left(\Cyl\right)}&=\left(\int^1_0\int^{2\pi}_0h(1+h\eta\kappa)|u(\eta,\theta)|^2\,\xdif\eta\,\xdif\theta\right)^{1/2},\\
 &=\|u\|_{\leb[2]{\mathcal{O}_h}},
 \intertext{and the $L^2$ norm of its exterior derivative $\xdif u$, denoted by $\|\xdif u\|_{\Lambda^1L^2_m}$ is equal to}
 \|\xdif u\|_{\Lambda^1L^2_\gfrak\left(\Cyl\right)}&=\left(\int^1_0\int^{2\pi}_0\frac{1+h\eta\kappa}{h}\left|\partial_\eta u(\eta,\theta)\right|^2+\frac{h}{1+h\eta\kappa}\left|\partial_\theta u(\eta,\theta)\right|^2\,\xdif\eta
 \,\xdif\theta\right)^{1/2},\\
 &=\|\grad u\|_{\leb[2]{\mathcal{O}_h}}.
 \end{align*}
To simplify our notations, for a 0-form $u$ defined on $\Cyl$, we define by $\|u\|_{H^1_\gfrak\left(\Cyl\right)}$ the following quantity
$$\|u\|_{H^1_\gfrak\left(\Cyl\right)}= \|u\|_{\Lambda^0L^2_m\left(\Cyl\right)}+ \|\xdif u\|_{\Lambda^1L^2_m\left(\Cyl\right)},$$
when the above integrals are well-defined. Observe that for a function $u\in\sob[1]{\mathcal{O}_h}$, we have:
$$\|u\|_{\sob[1]{\mathcal{O}_h}}=\|u\compo\Phi\|_{H^1_\gfrak\left(\Cyl\right)}.$$
\end{nota} 
\begin{rmrk}[Poincar\'e inequality in the thin layer]
Let $z$ belong to $H^1_\gfrak\left(\Cyl\right)$, such that \begin{align}\int^{2\pi}_0z(0,\theta)\,\xdif\theta=0.\label{gauge}\end{align}
Then, there exists an $h$-independant constant $C_{\mathcal{O}}$ such that
\begin{align}
\|z\|_{\Lambda^0L^2_\gfrak\left(\Cyl\right)}\leq C_{\mathcal{O}}\left\|\xdif z\right\|_{\Lambda^1L^2_\gfrak\left(\Cyl\right)}.\label{***}
\end{align}
We prove \eqref{***} using Fourier analysis. According to the definition  \eqref{h0} of $h_0$  there exists two constants $C_{\mathcal{O}}$ and $c_{\mathcal{O}}$ 
 depending on the domain $\mathcal{O}$ such that the following inequalities hold:
 \begin{subequations}
 \begin{align}
 \| z\|^2_{\Lambda^0L^2_\gfrak\left(D\right)}&\leq C_{\mathcal{O}}h\int^1_0\int^{2\pi}_0\left|z(\eta,\theta)\right|^2\,\xdif \theta\,\xdif\eta,\\
 \| \xdif z\|^2_{\Lambda^1L^2_\gfrak\left(D\right)}&\geq c_{\mathcal{O}}\left(\int^1_0\int^{2\pi}_0\frac{\left|\partial_\eta z(\eta,\theta)\right|^2}{h}+h\left|\partial_\theta z\right|^2\,\xdif \theta\,\xdif\eta\right).
 \end{align}\label{vi}
 \end{subequations}
For $k\in\Zr$, we denote by $\widehat{z}_k$  the $k^{\text{th}}$-Fourier coefficient (with respect to $\theta$) of $z$:
 \begin{align*}
 \widehat{z}_k=\frac{1}{2\pi}\int^{2\pi}_0z(\theta)\,e^{-\ri k\theta}\,\xdif\theta.
 \end{align*} 
 Since $\left(\widehat{\partial_\theta z}\right)_k=\ri k \widehat{z}_k $, we infer:
 \begin{align*}
 \forall k\neq0,\quad  \int^1_0\left|\widehat{z}_k(\eta)\right|^2\,\xdif\eta\leq \int^1_0\left|\left(\widehat{\partial_\theta z}\right)_k(\eta)\right|^2\,\xdif\eta .
 \end{align*}
 According to gauge condition \eqref{gauge}, we have:
  \begin{align*}
\widehat{z}_0(0)&=0,
  \intertext{thus, using the equality}
  \widehat{z}_0(\eta)&=\int^\eta_0\left(\widehat{\partial_\eta z}\right)_0(s)\xdif s,
  \intertext{we infer}
  \int^1_0\left|\widehat{z}_0(\eta)\right|^2\,\xdif\eta&\leq \int^1_0\left|\left(\widehat{\partial_\eta z}\right)_0(\eta)\right|^2\,\xdif\eta . 
  \intertext{Therefore,}
  \sum_{k\in\Zr}\int^1_0\left|\widehat{z}_k(\eta,\theta)\right|^2\,\xdif\eta&\leq \sum_{k\in\Zr}\left\{\int^1_0\left|\left(\widehat{\partial_\theta z}\right)_k(\eta)\right|^2\,\xdif\eta 
 +\int^1_0\left|\left(\widehat{\partial_\eta z}\right)_k(\eta)\right|^2 \right\}.
  \end{align*}
 We end the proof of \eqref{***} by using Parseval inequality and inequalities \eqref{vi}.
\end{rmrk}

\section{Asymptotic expansion of the steady state potential for an insulating membrane}\label{insulating membrane}
We derive asymptotic expansions with respect to $h$ of the potentials
$\left(V^c,V^m\right)$ solution to Problem \eqref{dielmetr}. The
membrane is insulating since the modulus of $\alpha$ is supposed
to be smaller than 1. However, our results are still valid if
$|\alpha|$ is bounded by a constant $C_0$ greater than 1. We emphasize
that the following results are valid for $\alpha$ tending to zero.

\subsection{Formal asymptotic expansion}\label{formal}
We write the following ansatz:
\begin{subequations}
\begin{align}
V^c&=V^c_0+hV^c_1+h^2V^c_2+\cdots,\label{ansatzVc}\\
V^m&=V^m_0+h V^m_1+h^2V^m_2+\cdots.\label{ansatzVm}
\end{align}\label{ansatz}
\end{subequations}
We multiply \eqref{laplacemetr} by $h^2(1+h\eta\kappa)^2$ and we order the powers of $h$ to obtain:
\begin{align}
&\forall (\eta,\theta)\in[0,1]\times\Tr,\notag\\
\begin{split}\partial^2_\eta V^m&+h\kappa\left\{3\eta\partial^2_\eta V^m+\partial_\eta V^m\right\}
+h^2\left\{3\eta^2\kappa^2\partial^2_\eta V^m+2\eta\kappa^2\partial_\eta V^m+\partial^2_\theta V^m\right\}\\
&+h^3\left\{\eta^3\kappa^3\partial^2_\eta V^m+\eta^2\kappa^3\partial_\eta V^m+\eta\kappa\partial^2_\theta V^m-\eta\kappa'\partial_\theta V^m\right\}=0
\end{split}\label{laplaceh}
\end{align}
We are now  ready to derive formally the terms of the asymptotic expansions of $V^c$ and $V^m$ by identifying the terms of the same power in $h$.

Recall that for $(m,n)$ in $\N^2$, $\delta_{m,n}$ is Kronecker symbol equal to 1 if $m=n$ and to 0 if $m\neq n$.
 By identifying the powers of $h$,
 we infer that for $l\in\N$, $V^c_l$ and $V^m_l$ satisfy the following equations:
\begin{subequations}
\begin{align}
\Delta V^c_l=&0,\text{  in $\mathcal{O}$},\label{pdeVcl}
\intertext{for all $(\eta,\theta)\in \Cyl$,}
\begin{split}\partial^2_\eta V^m_l=&-\Biggl\{\kappa\left\{3\eta\partial^2_\eta V^m_{l-1}+\partial_\eta V^m_{l-1}\right\}
\\&+3\eta^2\kappa^2\partial^2_\eta V^m_{l-2}+2\eta\kappa^2\partial_\eta V^m_{l-2}+\partial^2_\theta V^m_{l-2}
\\&+\eta^3\kappa^3\partial^2_\eta V^m_{l-3}+\eta^2\kappa^3\partial_\eta V^m_{l-3}+\eta\kappa\partial^2_\theta V^m_{l-3}
-\eta\kappa'\partial_\theta V^m_{l-3}\Biggr\},\end{split}\label{pdeVml}
\intertext{with transmission conditions}
\partial_n V^c_l\compo\Phi_0=&\alpha\left.\partial_\eta V^m_{l+1}\right|_{\eta=0},\label{CTNl}\\
V^c_l\compo\Phi_0=&\left.V^m_l\right|_{\eta=0},\label{CTDl}
\intertext{with boundary condition}
\left.\partial_\eta V^m_l\right|_{\eta=1}=&\delta_{l,1}f,\label{CLl}
\intertext{and with gauge condition}
\int_{\partial\mathcal{O}}V^c_l\xdif\sigma=&0.
\end{align}\label{dblaspt}
\end{subequations}
 In equations \eqref{dblaspt}, we have implicitly imposed 
\begin{equation}\left\{
\begin{aligned}
&V^c_l=0,\text{ if $l\leq -1$},\\
&V^m_l=0,\text{ if $l\leq -1$}.
\end{aligned}\right.\label{convention} 
\end{equation}
The next lemma ensures that for each non null integer $N$, the
functions $V^c_N$ and $V^m_N$ are entirely determined if the boundary
condition $\phi$ is enough regular. 
\begin{nota}
For $s\in\Er$, we denote by $\Cscr^\infty\left([0,1];\sob[s]{\Tr}\right)$ the space of functions $u$ defined for $(\eta,\theta)\in[0,1]\times\Tr$, such that for almost all $\theta\in\Tr$, $u(\cdot,\theta)$
 belongs to $\Cscr^{\infty}\left([0,1]\right)$, and such that for all $\eta\in[0,1]$, $u(\eta,\cdot)$ belongs to $\sob[s]{\Tr}$.
\end{nota}

\begin{lemm}\label{lemme1}
We suppose that $\partial \mathcal{O}$ is smooth.

For $N\in\N$ and $p\geq 0$ we suppose that $\phi$ belongs to $\sob[N+p-1/2]{\partial\W_h}$ and let $|\alpha |\leq1$.

Then the functions $V^m_0,\cdots,V^m_N$ and $V^c_0,\cdots,V^c_{N}$ are uniquely determined and they belong to the respective functional spaces:
\begin{subequations}
\begin{align}
\forall k=0,\cdots, N,\notag\\
&V^m_k\in\Cscr^\infty\left([0,1];\sob[N+p-k+1/2]{\Tr}\right),\\
&V^c_{k}\in\sob[N+p-k+1]{\mathcal{O}}.
\end{align}\label{regularite}
\end{subequations}
Moreover, there exists a constant $C_{N,\mathcal{O},p}$ such that:
\begin{subequations}
\begin{align}
\forall k=0,\cdots, N,\notag\\
\sup_{\eta\in[0,1]}\left\|V^m_k(\eta,\cdot)\right\|_{\sob[N+p-k+1/2]{\Tr}}&\leq C_{N,\mathcal{O},p}\|\frm\|_{\sob[N+p-1/2]{\partial\mathcal{O}}},\\
\left\|V^c_k\right\|_{\sob[N+p-k+1]{\mathcal{O}}}&\leq |\alpha|C_{N,\mathcal{O},p}\|\frm\|_{\sob[N+p-1/2]{\partial\mathcal{O}}}.
\end{align}\label{estimVmkVck}
\end{subequations}
 \end{lemm}
 \begin{rmrk}
 To simplify, we suppose that $|\alpha|\leq 1$, but the same result may be obtained if there exists $C_0>1$ such that $|\alpha|\leq C_0$. In this case, the constant $C_{N,\mathcal{O},p}$ would also 
 depends on $C_0$.
 \end{rmrk}

\begin{proof}
Since $\partial\mathcal{O}$ is smooth and since $\phi$ belongs to $\sob[N+p-1/2]{\partial\W_h}$, for $N\geq0$ and $p\geq0$,  then the functions $f$ and $\frm$ defined by \eqref{ftangentielle} and by 
\eqref{frm} belong respectively to $\sob[N+p-1/2]{\Tr}$
and to $\sob[N+p-1/2]{\partial\mathcal{O}}$.
We prove this lemma by recursive process.

\begin{itemize}
\item $N=0$. Let $p\geq 0$ and let $\phi$ belong to $\sob[p-1/2]{\partial\W_h}$. 
\end{itemize}
Thus $f$ and $\frm$ belong respectively to $\sob[p-1/2]{\Tr}$ and
$\sob[p-1/2]{\partial\mathcal{O}}$. Using \eqref{pdeVml} and \eqref{CLl}, we infer:
\begin{align}
\begin{cases}\left.\partial^2_\eta V^m_0=0\right.,\\
\left.\partial_\eta V^m_0\right|_{\eta=1}=0,\end{cases}\label{Vm01deriv}
\end{align}
hence, $\partial_\eta V^m_0= 0$. According to \eqref{pdeVml} and to \eqref{CLl},  we straight infer 
$$\partial_\eta V^m_1=f.$$
Therefore by \eqref{pdeVcl} and \eqref{CTNl} the function $V^c_0$ satisfies the following Laplace problem:
\begin{subequations}
\begin{align}
&\left.\Delta V^c_0\right.=0,\\
&\left.\partial_n V^c_0\right|_{\partial\mathcal{O}}=\alpha \frm,
\intertext{with gauge condition}
&\int_{\partial\mathcal{O}}V^c_0\xdif \sigma=0.
\end{align}\label{Vc0}
\end{subequations}
According to \eqref{CTDl}, we infer \begin{align}V^m_0=V^c_0\compo\Phi_0\label{Vm0},\end{align}
hence $V^c_0$ and $V^m_0$ are entirely determined and they belong to the following spaces:
\begin{align*}
&V^m_0\in\Cscr^\infty\left([0,1];\sob[p+1/2]{\Tr}\right),\\
&V^c_0\in\sob[p+1]{\mathcal{O}}.
\end{align*}
Observe also that there exists a constant $C_{\mathcal{O},p}$ such that 
\begin{align*}
\sup_{\eta\in[0,1]}\left\|V^m_0(\eta,\cdot)\right\|_{\sob[p+1/2]{\Tr}}&\leq C_{\mathcal{O},p}\|\frm\|_{\sob[p-1/2]{\partial\mathcal{O}}},\\
\left\|V^c_0\right\|_{\sob[p+1]{\mathcal{O}}}&\leq |\alpha|C_{\mathcal{O},p}\|\frm\|_{\sob[p-1/2]{\partial\mathcal{O}}}.
\end{align*}

\begin{itemize}
\item Induction.
\end{itemize} Let $N\geq0$. Suppose that for all $p\geq 0$, for all $\phi\in\sob[N+p-1/2]{\partial\W_h}$ and  for $M=0,\cdots, N$ the functions $V^c_M$ and $V^m_M$ are known.
Suppose that 
they belong respectively to 
$\sob[N+p-M+1]{\mathcal{O}}$ and to  $V^m_M\in\Cscr^\infty\left([0,1];\sob[N+p-M+1/2]{\Tr}\right)$ and that estimates \eqref{estimVmkVck} hold.

Let $\phi$ belong to $\sob[N+p+1/2]{\partial\W_h}$. Therefore, for $M=0,\cdots,N,$ the functions $V^c_M$ and $V^m_M$ are known, they belong respectively to 
$\sob[N+p-M+2]{\mathcal{O}}$ and to  $V^m_M\in\Cscr^\infty\left([0,1];\sob[N+p-M+3/2]{\Tr}\right)$ and the following estimates hold:
\begin{align*}
\forall M=0,\cdots, N,\notag\\
\sup_{\eta\in[0,1]}\left\|V^m_M(\eta,\cdot)\right\|_{\sob[N+p-M+3/2]{\Tr}}&\leq C_{N,\mathcal{O},p}\|\frm\|_{\sob[N+p+1/2]{\partial\mathcal{O}}},\\
\left\|V^c_M\right\|_{\sob[N+p-M+2]{\mathcal{O}}}&\leq |\alpha|C_{N,\mathcal{O},p}\|\frm\|_{\sob[N+p+1/2]{\partial\mathcal{O}}}.
\end{align*}
We are going to build $V^c_{N+1}$ and $V^m_{N+1}$.
From \eqref{pdeVml} and \eqref{CLl}, we infer, for all $(\eta,\theta)\in \Cyl$,
\begin{align*}
\begin{split}\partial^2_\eta V^m_{N+1}=&-\Biggl\{\kappa\left\{3\eta\partial^2_\eta V^m_{N}+\partial_\eta V^m_{N}\right\}
\\&+3\eta^2\kappa^2\partial^2_\eta V^m_{N-1}+2\eta\kappa^2\partial_\eta V^m_{N-1}+\partial^2_\theta V^m_{N-1}
\\&+\eta^3\kappa^3\partial^2_\eta V^m_{N-2}+\eta^2\kappa^3\partial_\eta V^m_{N-2}+\eta\kappa\partial^2_\theta V^m_{N-2}
-\eta\kappa'\partial_\theta V^m_{N-2}\Biggr\},\end{split}
\\
\left.\partial_\eta V^m_{N+1}\right|_{\eta=1}&=0.
\end{align*}
Recall that we use convention \eqref{convention}. Since we have supposed that $V^m_{M}$ is known for $M\leq N$ and belongs to $\Cscr^\infty\left([0,1];\sob[N+1+p-M-1/2]{\Tr}\right)$, 
we infer that:
\begin{align}
&\forall (s,\theta)\in\Cyl,\notag\\
\begin{split}\partial_\eta V^m_{N+1}(s,\cdot)=&\int^1_s\biggl\{\kappa\left\{3\eta\partial^2_\eta V^m_{N}+\partial_\eta V^m_{N}\right\}
\\&+3\eta^2\kappa^2\partial^2_\eta V^m_{N-1}+2\eta\kappa^2\partial_\eta V^m_{N-1}+\partial^2_\theta V^m_{N-1}
\\&+\eta^3\kappa^3\partial^2_\eta V^m_{N-2}+\eta^2\kappa^3\partial_\eta V^m_{N-2}+\eta\kappa\partial^2_\theta V^m_{N-2}
-\eta\kappa'\partial_\theta V^m_{N-2}\biggr\}\xdif \eta,\end{split}\label{partialetaVmN+1}
\end{align}
is entirely determined and belongs to $\Cscr^\infty\left([0,1];\sob[p+1/2]{\Tr}\right)$. Moreover, since $\partial_\eta V^m_{N+1}$ is known, we infer exactly by the same way that $\partial_\eta V^m_{N+2}$ is also determined. Actually, it is equal to
\begin{align*}
&\forall (s,\theta)\in\Cyl,\\
\begin{split}\partial_\eta V^m_{N+2}(s,\cdot)=&\int^1_s\biggl\{\kappa\left\{3\eta\partial^2_\eta V^m_{N+1}+\partial_\eta V^m_{N+1}\right\}
\\&+3\eta^2\kappa^2\partial^2_\eta V^m_{N}+2\eta\kappa^2\partial_\eta V^m_{N}+\partial^2_\theta V^m_{N}
\\&+\eta^3\kappa^3\partial^2_\eta V^m_{N-1}+\eta^2\kappa^3\partial_\eta V^m_{N-1}+\eta\kappa\partial^2_\theta V^m_{N-1}
-\eta\kappa'\partial_\theta V^m_{N-1}\biggr\}\xdif \eta,\end{split}
\end{align*}
and it belongs to $\Cscr^\infty\left([0,1];\sob[p+1/2]{\Tr}\right)$.
According to \eqref{CTNl}, the function $V^c_{N+1}$ is then uniquely determined by 
\begin{subequations}
\begin{align}
&\left.\Delta V^c_{N+1}\right.=0,\\
&\left.\partial_n V^c_{N+1}\right|_{\partial\mathcal{O}}=\alpha \partial_\eta V^m_{N+2}\compo\Phi^{-1}_0,
\intertext{with gauge condition}
&\int_{\partial\mathcal{O}}V^c_{N+1}\xdif \sigma=0.
\end{align}\label{Vcl+1}
\end{subequations}
Moreover, it belongs to $\sob[p+1]{\mathcal{O}}$.
Transmission condition \eqref{CTDl} implies the following expression of $V^m_{N+1}$:
$$\forall s\in(0,1),\quad V^m_{N+1}(s,\cdot)=\int^s_0 \partial_\eta V^m_{N+1}(\eta,\cdot)\,\xdif\eta+V^c_{N+1}\compo\Phi_0,$$
where $\partial_\eta V^m_{N+1}$ is given by \eqref{partialetaVmN+1} and belongs to $\Cscr^\infty\left([0,1];\sob[p+1/2]{\Tr}\right)$. We infer also that there exists $C_{N+1,\mathcal{O},p}>0$ such that
\begin{align*}
\sup_{\eta\in[0,1]}\left\|V^m_{N+1}(\eta,\cdot)\right\|_{\sob[p+1/2]{\Tr}}&\leq C_{N+1,\mathcal{O},p}\|\frm\|_{\sob[N+p+1/2]{\partial\mathcal{O}}},\\
\left\|V^c_{N+1}\right\|_{\sob[p+1]{\mathcal{O}}}&\leq |\alpha|C_{N+1,\mathcal{O},p}\|\frm\|_{\sob[N+p+1/2]{\partial\mathcal{O}}},
\end{align*}
 hence the lemma.
\end{proof}
Observe that the functions $(V^c_k,V^m_k)$ are these given in Theorem~\ref{thm1}. 
 \subsection{Error Estimates of Theorem~\ref{thm1}}\label{estimates}
Let us prove now the estimates of Theorem~\ref{thm1}. 
Let $N\in\N$ and $\phi$ belong to $\sob[N+3/2]{\partial \W_h}$. The
 function $\frm$ is defined by \eqref{frm}.
Let $R^c_N$ and $R^m_N$ be the functions defined by:
 \begin{align*}
\begin{cases} R^c_N=V_h-\sum^N_{k=0}V^c_kh^k,\,\text{in $\mathcal{O}$},\\
 R^m_N=V_h\compo\Phi-\sum^N_{k=0}V^m_kh^k,\,\text{in $\Cyl$}.
 \end{cases}
 \end{align*}
We have to prove that there exists a constant $C_{\mathcal{O},N}>0$ depending only on the domain $\mathcal{O}$ and on $N$ such that 
\begin{subequations}
 \begin{align}
\|R^c_N\|_{\sob[1]{\mathcal{O}}}&\leq C_{\mathcal{O},N}\|\frm\|_{\sob[N+3/2]{\partial\mathcal{O}}}|\alpha| h^{N+1/2},\\
\|R^m_N\|_{H^1_\gfrak\left(\Cyl\right)}&\leq C_{\mathcal{O},N}\|\frm\|_{\sob[N+3/2]{\partial\mathcal{O}}}h^{N+1/2}.
 \end{align}\label{resultat1}
\end{subequations}
Moreover, if $\phi$ belongs to $\sob[N+5/2]{\partial \W_h}$, then we have
\begin{subequations}
 \begin{align}
\|R^c_N\|_{\sob[1]{\mathcal{O}}}&\leq C_{\mathcal{O},N}\|\frm\|_{\sob[N+5/2]{\partial\mathcal{O}}}|\alpha| h^{N+1},\\
\|R^m_N\|_{H^1_\gfrak\left(\Cyl\right)}&\leq C_{\mathcal{O},N}\|\frm\|_{\sob[N+5/2]{\partial\mathcal{O}}}h^{N+1/2}.
 \end{align}\label{resultat2}
\end{subequations}
\begin{proof}[Proof of Theorem~\ref{thm1}]
Since $\phi$ belongs to $\sob[N+3/2]{\partial \W_h}$, according to the previous lemma, the couples of  functions $\left(R^c_{N},R^m_{N}\right)$ and  $\left(R^c_{N+1},R^m_{N+1}\right)$
 are well defined and belong to 
$\sob[1]{\mathcal{O}}\times H^1_\gfrak\left(\Cyl\right)$. The Sobolev
space $H^1_\gfrak\left(\Cyl\right)$ is defined in Notation
\ref{notasob}.

Denote by $g_N$ the following function defined on $\Cyl$:
\begin{align}
\begin{split}g_N=&\kappa\left(3\eta\partial^2_\eta V^m_{N}+\partial_\eta V^m_{N}\right)
+3\eta^2\kappa^2\partial^2_\eta V^m_{N-1}+2\eta\kappa^2\partial_\eta V^m_{N-1}+\partial^2_\theta V^m_{N-1}\\&
+\eta^3\kappa^3\partial^2_\eta V^m_{N-2}+\eta^2\kappa^3\partial_\eta V^m_{N-2}+\eta\kappa\partial^2_\theta V^m_{N-2}
-\eta\kappa'\partial_\theta V^m_{N-2}\\&+
h\Biggl(3\eta^2\kappa^2\partial^2_\eta V^m_{N}+2\eta\kappa^2\partial_\eta V^m_{N}+\partial^2_\theta V^m_{N}
\\&+\eta^3\kappa^3\partial^2_\eta V^m_{N-1}+\eta^2\kappa^3\partial_\eta V^m_{N-1}+\eta\kappa\partial^2_\theta V^m_{N-1}
-\eta\kappa'\partial_\theta V^m_{N-1}
\Biggr)\\&+
h^2\Biggl(\eta^3\kappa^3\partial^2_\eta V^m_{N}+\eta^2\kappa^3\partial_\eta V^m_{N}+\eta\kappa\partial^2_\theta V^m_{N}
-\eta\kappa'\partial_\theta V^m_{N}
\Biggr)\end{split}
\end{align}
According to the previous lemma and since $\phi$ belongs to $\sob[N+1/2]{\partial\W_h}$, the above function $g_N$ belongs to 
$\Cscr^\infty\left([0,1];\sob[-1/2]{\Tr}\right)$ and the function $\partial_\eta V^m_N$ belongs to 
$\Cscr^\infty\left([0,1];\sob[3/2]{\Tr}\right)$.  Moreover, there exists a constant $C_{N,\mathcal{O}}$ such that
\begin{align}
\begin{cases}
\sup_{\eta\in[0,1]}\|g_N(\eta,\cdot)\|_{H^{-1/2}\left(\Tr\right)}\leq C_{N,\mathcal{O}}\|f\|_{\sob[N+1/2]{\Tr}},\\
\sup_{\eta\in[0,1]}\left\|\partial_\eta V^m_N(\eta,\cdot)\right\|_{H^{3/2}\left(\Tr\right)}\leq C_{N,\mathcal{O}}\|f\|_{\sob[N+1/2]{\Tr}}.\end{cases}
\label{estimgN}
\end{align}
The functions $R^c_N$ and $R^m_N$ satisfy the following problem:
\begin{align*}
&\Delta R^c_N=0,\text{ in $\mathcal{O}$},\\
&\partial_\eta\left(\frac{1+h\eta\kappa}{h}\partial_\eta R^m_N\right)+\partial_\theta\left(\frac{h}{1+h\eta\kappa}\partial_\theta R^m_N\right)=\frac{-h^N}{(1+h\eta\kappa)}g_N,
\intertext{with transmission conditions:}
&\partial_nR^c_N\compo\Phi_0=\frac{\alpha}{h}\left(\left.\partial_\eta R^m_N\right|_{\eta=0}+h^{N+1}\left.\partial_\eta V^m_N\right|_{\eta=0}\right),\\
&R^c_N\compo\Phi_0=\left.R^m_N\right|_{\eta=0},
\intertext{with boundary condition}
&\left.\partial_\eta R^m_N\right|_{\eta=1}=0,
\intertext{and with gauge condition}
&\int_{\partial\mathcal{O}}R^c_N\,\xdif\sigma=0.
\end{align*}
By multiplying the above equality by $\overline{R_N}$ and by integration by parts, we infer that:
\begin{align}\begin{split}\left\|\xdif R^c_N\right\|^2_{\Lambda^1\leb[2]{\mathcal{O}}}+\alpha\left\|\xdif R^m_N\right\|^2_{\Lambda^1L^2_\gfrak\left(\Cyl\right)}&
=-\alpha h^N\int_{\Cyl}g_N(\eta,\theta)\overline{R^m_N}(\eta,\theta)\,\xdif\eta\,\xdif\theta\\
&+\alpha h^{N+1}\int_{\Tr}\left.\partial_\eta V^m_N\right|_{\eta=0}\left.\overline{R^m_N}\right|_{\eta=0}\,\xdif\theta\\
&-\alpha h^{N+1}\int_{\Tr}\kappa\left.\partial_\eta V^m_N\right|_{\eta=1}\left.\overline{R^m_N}\right|_{\eta=1}\,\xdif\theta.\end{split}\label{equality}
 \end{align} 
By hypothesis \eqref{hyp}, and using by Cauchy-Schwarz inequality and estimates \eqref{estimgN}, we infer that there exists a constant $C_{\mathcal{O},N}>0$ such that  
\begin{align*}\Re(\alpha)\left\|\xdif R^m_N\right\|^2_{\Lambda^1L^2_\gfrak\left(\Cyl\right)}&\leq
|\alpha|C_{\mathcal{O},N}h^{N-1/2}\|R^m_N\|_{H^1_\gfrak\left(\Cyl\right)}\|f\|_{\sob[N+1/2]{\Tr}}\|R^m_N\|_{H^1_\gfrak\left(\Cyl\right)},
\intertext{and}
\left|\Im(\alpha)\right|\left\|\xdif R^m_N\right\|^2_{\Lambda^1L^2_\gfrak\left(\Cyl\right)}&\leq
|\alpha|C_{\mathcal{O},N}h^{N-1/2}\|R^m_N\|_{H^1_\gfrak\left(\Cyl\right)}\|f\|_{\sob[N+1/2]{\Tr}}\|R^m_N\|_{H^1_\gfrak\left(\Cyl\right)},
\intertext{hence}
\left\|\xdif R^m_N\right\|^2_{\Lambda^1L^2_\gfrak\left(\Cyl\right)}&\leq C_{\mathcal{O},N}h^{N-1/2}\|f\|_{\sob[N+1/2]{\Tr}}\|R^m_N\|_{H^1_\gfrak\left(\Cyl\right)}.
 \end{align*} 
Since $\int_{\Tr}\left.R^m_N\right|_{\eta=0}\xdif\theta=0$, by Poincar\'e inequality \eqref{***}, there exists a strictly positive constant $C_\mathcal{O}$, which does not depend on $h$ such 
that 
$$\|R^m_N\|_{\Lambda^0L^2_\gfrak\left(\Cyl\right)}\leq C_{\mathcal{O}}\left\|\xdif R^m_N\right\|_{\Lambda^1L^2_\gfrak\left(\Cyl\right)},$$
hence \begin{align*}\left\| R^m_N\right\|_{H^1_\gfrak\left(\Cyl\right)}\leq
C_{\mathcal{O},N}\|f\|_{\sob[N+1/2]{\Tr}}h^{N-1/2},
\intertext{and therefore we deduce directly from the above estimate and from \eqref{equality},}
\left\| R^c_N\right\|_{\sob[1]{\mathcal{O}}}\leq
C_{\mathcal{O},N}\|f\|_{\sob[N+1/2]{\Tr}}|\alpha|h^{N-1/2}.
 \end{align*} 
The above estimate holds for $\phi\in\sob[N+1/2]{\partial\W_h}$. Since $\phi$ belongs to $\sob[N+3/2]{\partial\W_h}$, we obtain the same result by replacing $N$ by $N+1$:
\begin{equation}\left\{\begin{aligned}
&\left\| R^m_{N+1}\right\|_{H^1_\gfrak\left(\Cyl\right)}\leq
C_{\mathcal{O},N+1}\|f\|_{\sob[N+3/2]{\Tr}}h^{N+1/2},\\
&\left\| R^c_{N+1}\right\|_{\sob[1]{\mathcal{O}}}\leq
C_{\mathcal{O},N+1}\|f\|_{\sob[N+3/2]{\Tr}}|\alpha|h^{N+1/2}.
 \end{aligned}\right.
\end{equation} 
According to the previous lemma, the functions $V^c_{N+1}$ and $V^m_{N+1}$ are well-defined and there exists a constant $C_{N,\mathcal{O}}$ such that:
 \begin{align*}
 &\|V^c_{N+1}\|_{\sob[1]{\mathcal{O}}}\leq |\alpha|C_{N,\mathcal{O}}\|\frm\|_{\sob[N+3/2]{\partial\mathcal{O}}},\\
 &\|V^m_{N+1}\|_{H^1_\gfrak\left(\Cyl\right)}\leq \frac{C_{N,\mathcal{O}}}{\sqrt{h}}\|\frm\|_{\sob[N+3/2]{\partial\mathcal{O}}}.
 \end{align*}
Writing $$R^c_N=R^c_{N+1}+V^c_{N+1}h^{N+1},$$ 
and $$R^m_N=R^m_{N+1}+V^m_{N+1}h^{N+1},$$
 we infer that 
$$\left\| R^c_N\right\|_{\sob[1]{\mathcal{O}}}\leq
C_{\mathcal{O},N}\|f\|_{\sob[N+3/2]{\Tr}}|\alpha|h^{N+1/2},$$
and 
$$\left\| R^m_N\right\|_{H^1_{\gfrak}\left(\Cyl\right)}\leq
C_{\mathcal{O},N}\|f\|_{\sob[N+3/2]{\Tr}}h^{N+1/2}.$$
If $\phi$ belongs to $\sob[N+5/2]{\partial\W_h}$, we write $$R^c_N=R^c_{N+2}+V^c_{N+1}h^{N+1}+V^c_{N+2}h^{N+2},$$ 
and $$R^m_N=R^m_{N+1}+V^m_{N+1}h^{N+1}+V^m_{N+2}h^{N+2},$$ to obtain estimates \eqref{resultat2},
hence Theorem~\ref{thm1}.
\end{proof}
\begin{rmrk}[The case of an insulating inner domain]
Consider Problem~\eqref{diel11}:
  \begin{align*}
   \dverg\left(\gamma_h\grad u\right)&=0
  \text{ in $\W_h$},\\
   \frac{\partial u}{\partial n}&=\phi
  \text{ on $\partial\W_h$},\\
  \int_{\partial \mathcal{O}}u\,\xdif \sigma&=0.
  \end{align*}
 If the inner domain is perfectly insulating (\textit{i.e.} if
 $\gamma_h$ vanishes in $\mathcal{O}$), the steady state potential in the membrane satisfies:
\begin{align*}
&\frac{1}{h^2}\partial_\eta\left((1+h\eta\kappa)\partial_\eta u^m\right)+\partial_\theta\left(\frac{1}{1+h\eta\kappa}\partial_\theta u^m\right)=0,\text{ in $\Cyl$},\\
\intertext{with the following boundary conditions:}
&\left.\partial_\eta u^m\right|_{\eta=0}=0,\quad \left.\partial_\eta u^m\right|_{\eta=1}=hf,
\intertext{and with gauge condition}
&\int^{2\pi}_0u^m|_{\eta=0}\xdif\theta=0.
\end{align*} 
By identifying the terms of the same power of $h$ we would obtain:
$u^m_0=0$, and $u^m_1$ would satisfy:
\begin{align*}
&\partial^2_\eta u^m_1=0,\text{ in $\Cyl$,}\\
&\left.\partial_\eta u^m_1\right|_{\eta=0}=0,\,\left.\partial_\eta u^m_1\right|_{\eta=1}=f,\\
&\int^{2\pi}_0u^m_1|_{\eta=0}\xdif\theta=0,
\end{align*} 
which is a non sense as soon as $f\neq0$. Our ansatz \eqref{ansatz} fails.
Actually, the asymptotic
expansion of $u^m$ begins at the order $-1$: a boundary layer
phenomenon appears.This is described in the next section. 
\end{rmrk}

\section{Asymptotic expansion of the steady state potential for an insulating inner domain}\label{insulating inner domain}
Consider now the solution $u_h$ to Problem \eqref{diel11}.
In this section, we suppose  
\begin{subequations}
 \begin{align}
&|\alpha|\text{ tends to zero},\\
&\Re(\alpha)> 0\,\text{ or }\,\Bigl\{\Re(\alpha)=0\,\text{ and }\,\Im(\alpha)\neq0\Bigr\}.
 \end{align}\label{alpha}
\end{subequations}
Thus the inner domain is insulating.
Let $\beta$ be a complex parameter satisfying:
$$Re(\beta)>0, or \left(Re(\beta)=0,\text{ and } \Im(\beta)\neq0\right).$$
The modulus of $\beta$ may tend to infinity, or to zero but it must satisfy:
$$|\beta|=o\left(\frac{1}{h}\right),\quad\text{and}\quad \frac{1}{|\beta|}=o\left(\frac{1}{h}\right).$$

We suppose that $u$ may be written as follows:
$$u_h=\frac{1}{h}u_{-1}+u_0+hu_1+\cdots.$$
We denote by $u^c$ and $u^m\compo\Phi^{-1}$ the respective restrictions
of $u_h$ to $\mathcal{O}$ and to $\mathcal{O}_h$.
One of the two following cases holds.
\begin{hypo}[$\alpha=\beta h^q$]\label{hypo1}
There exists $q\geq 1$ such that:
\begin{align}
\alpha&=\beta h^q.
\end{align} 
\end{hypo}
\begin{hypo}[$\alpha=o(h^{N})$, $\forall N\in\N$]\label{hypo2}
The complex parameter $\alpha$ satisfies \eqref{alpha} and for all $N\in\N$,
\begin{align}
\forall N\in\N,\quad |\alpha|=o( h^N).
\end{align}
\end{hypo}
First we suppose that Hypothesis~\ref{hypo1} holds: we will discuss on Hypothesis~\ref{hypo2} later on.
We denote by $(u^{c,q},u^{m,q})$ the solution to Problem~\ref{diel11}
under the Hypothesis~\ref{hypo1}.

According to \eqref{dielmetru}, by ordering and identifying the terms
of the same power of $h$, for $k\in\N\cup{-1}$, for $q\in\N^*$, $u^{c,q}_k$ and $u^{m,q}_k$ satisfy:
\begin{subequations}
\begin{align}
\Delta u^{c,q}_l=&0,\text{  in $\mathcal{O}$},\label{pdeucl}
\intertext{for all $(\eta,\theta)\in \Cyl$,}
\begin{split}\partial^2_\eta u^{m,q}_l=&-\Biggl\{\kappa\left\{3\eta\partial^2_\eta u^{m,q}_{l-1}+\partial_\eta u^{m,q}_{l-1}\right\}
\\&+3\eta^2\kappa^2\partial^2_\eta u^{m,q}_{l-2}+2\eta\kappa^2\partial_\eta u^{m,q}_{l-2}+\partial^2_\theta u^{m,q}_{l-2}
\\&+\eta^3\kappa^3\partial^2_\eta u^{m,q}_{l-3}+\eta^2\kappa^3\partial_\eta u^{m,q}_{l-3}+\eta\kappa\partial^2_\theta u^{m,q}_{l-3}
-\eta\kappa'\partial_\theta u^{m,q}_{l-3}\Biggr\},\end{split}\label{pdeuml}\\
u^{c,q}_l\compo\Phi_0=&\left.u^{m,q}_l\right|_{\eta=0},\label{CTDlcytins}\\
\left.\partial_\eta u^{m,q}_l\right|_{\eta=1}=&\delta_{l,1}f,\label{CLlcytins}\\
\int_{\partial\mathcal{O}}u^{c,q}_l\xdif\sigma=&0.
\intertext{Transmission condition \eqref{transmi1u} coupled
  with
  Hypothesis~\ref{hypo1} implies:}
\beta \partial_n u^{c,q}_{l-1-q}\compo\Phi_0=&\left.\partial_\eta
u^{m,q}_{l}\right|_{\eta=0},\label{CTNlcytins}
\end{align}\label{dblasptcytins}
\end{subequations}
In equations \eqref{dblasptcytins}, we have implicitly imposed 
\begin{align}\begin{cases}
u^{c,q}_l=0,\text{ if $l\leq -2$},\\
u^{m,q}_l=0,\text{ if $l\leq -2$}.\end{cases}\label{conventioncytins}
\end{align}

Let us now derive formal asymptotics of $u$ when Hypothesis~\ref{hypo1} holds.

\subsection{Formal asymptotics}
\begin{itemize}
\item $N=-1$.
\end{itemize}

The functions $u^{m,q}_{-1}$ satisfies
\begin{equation*}
\left\{\begin{aligned}
&\partial^2_\eta u^{m,q}_{-1}=0,\text{in $\Cyl$},\\
&\partial_\eta u^{m,q}_{-1}|_{\eta=0}=0,\quad \partial_\eta u^{m,q}_{-1}|_{\eta=1}=0,
\end{aligned}\right.
\end{equation*}
hence $u^{m,q}_{-1}$ depends only on the variable $\theta$. 
Observe that we have, for almost all $\theta\in\Tr$ the following equality:
$$u^{m,q}_{-1}(\theta)=u^{c,q}_{-1}\compo\Phi_0(\theta).$$
\begin{itemize}
\item $N=0$.
\end{itemize}
The function $u^{m,q}_0$ satisfies:
\begin{equation*}
\left\{\begin{aligned}
&\partial^2_\eta u^{m,q}_{0}=0,\text{in $\Cyl$},\\
&\partial_\eta u^{m,q}_{0}|_{\eta=0}=0,\quad \partial_\eta u^{m,q}_{0}|_{\eta=1}=0,
\end{aligned}\right.
\end{equation*}
hence, $\partial_\eta u^{m,q}_0$ vanishes identically in $\Cyl$. 
\begin{itemize}
\item $N=1$.
\end{itemize}
The functions $u^{m,q}_1$ satisfy:
\begin{equation*}
\left\{\begin{aligned}
&\partial^2_\eta u^{m,q}_{1}=-\partial^2_\theta u^{m,q}_{-1},\text{in $\Cyl$},\\
&\partial_\eta u^{m,q}_{1}|_{\eta=0}=\beta\partial_n u^{c,q}_{-1-q}\compo\Phi_0,\quad \partial_\eta u^{m,q}_{1}|_{\eta=1}=f.
\end{aligned}\right.
\end{equation*}
Therefore for $q=1$ we obtain the following equality:
$$-\partial^2_\theta u^{m,1}_{-1} +\beta\partial_n u^{c,1}_{-1}\compo\Phi_0=f,$$
hence the following boundary condition imposed to $u^{c,1}_{-1}$ on $\partial\mathcal{O}$:
$$-\left.\partial^2_\theta u^{c,1}_{-1}\right|_{\partial\mathcal{O}} +\beta\left.\partial_n u^{c,1}_{-1}\right|_{\partial\mathcal{O}}=\frm.$$
Therefore, the function $u^{c,1}_{-1}$ is solution to the following problem:
\begin{equation}
\left\{\begin{aligned}
&\Delta u^{c,1}_{-1}=0,\text{ in $\mathcal{O}$,}\\
&-\left.\partial^2_\theta u^{c,1}_{-1}\right|_{\partial\mathcal{O}} +\beta\left.\partial_n u^{c,1}_{-1}\right|_{\partial\mathcal{O}}=\frm,\\
&\int_{\partial\mathcal{O}}u^{c,1}_{-1}\xdif_{\partial\mathcal{O}}=0.
\end{aligned}\right.\label{uc-1},
\end{equation}
\begin{equation}
\forall (\eta,\theta)\in\Cyl,\quad u^{m,1}_{-1}=u^{c,1}_{-1}|_{\partial\mathcal{O}}\compo\Phi_0.\label{um-1}
\end{equation}
Since $\Re(\beta)>0$, a straight application of Lax-Milgram theorem ensures that $u^{c,1}$ is uniquely determined and belongs to 
$\sob[1]{\mathcal{O}}$ as soon as the boundary data belongs to $\sob[-3/2]{\partial\mathcal{O}}$.

If $q\geq 2$, the function $u^{m,q}_1$ satisfies:
\begin{align}
-\partial^2_\theta u^{m,q}_{-1}=f.\label{um2-1}\end{align}
Since $\int_\Tr u^{m,q}_{-1}\xdif\theta=0$, equality~\eqref{um2-1} defines
uniquely $u^{m,q}_{-1}$.
We infer that $u^{c,q}_{-1}$ is solution to the following problem:
\begin{equation}
\left\{\begin{aligned}
&\Delta u^{c,q}_{-1}=0,\text{ in $\mathcal{O}$,}\\
&-\left. u^{c,q}_{-1}\right|_{\partial\mathcal{O}}=u^{m,q}_{-1}\compo\Phi^{-1}_0 .
\end{aligned}\right.\label{uc2-1}
\end{equation}
Hence we have determined $u^{m,q}_{-1}$ and $u^{c,q}_{-1}$ for $q\in\N^*$. 
Observe that $u^{c,1}_{-1}$ is solution to Laplace equation with mixed
boundary condition, and for $q\geq2$ the potential 
$u^{c,q}_{-1}$ is the solution to Laplace 
equation with Dirichlet boundary condition, 
while for an insulating membrane, we obtained Neumann conditions for the approximated 
steady state potentials.

Let us now determined $u^{m,q}_{N}$ and $u^{c,q}_{N}$ 
for $q\in\N^*$ by recurrence.
\begin{itemize}
\item Induction.
\end{itemize}
Suppose that for $N\geq 0$, the functions $u^{m,q}_{N-1}$, $u^{c,q}_{N-1}$, $\partial_\eta u^{m,q}_{N}$ and $\partial_\eta u^{m,q}_{N+1}$
 are built. 

The function $u^{m,q}_{N+2}$ satisfies:
\begin{equation*}
\left\{\begin{aligned}
&\partial^2_\eta u^{m,q}_{N+2}=-\kappa\left(3\eta\partial^2_\eta u^{m,q}_{N+1}+\partial_\eta u^{m,q}_{N+1}\right)-\partial^2_\theta u^{m,q}_{N}-\eta\kappa\partial^2_\theta u^{m,q}_{N-1}+\eta\kappa'
\partial_\theta u^{m,q}_{N-1}
,\text{in $\Cyl$},\\
&\partial_\eta u^{m,q}_{N+2}|_{\eta=0}=\beta\partial_n u^{c,q}_{N+1-q}\compo\Phi_0,\quad \partial_\eta u^{m,q}_{N+2}|_{\eta=1}=0.
\end{aligned}\right.
\end{equation*}
Denote by $\phi^q_{N}$ the following function:
$$\phi^q_N=\int^1_0\Bigl(\kappa\left(3\eta\partial^2_\eta u^{m,q}_{N+1}+\partial_\eta u^{m,q}_{N+1}\right)+\eta\kappa\partial^2_\theta u^{m,q}_{N-1}-\eta\kappa'
\partial_\theta u^{m,q}_{-1}\Bigr)\,\xdif\eta.$$
Since $\partial^2_\eta u^{m,q}_{N+1}$ and $\partial_\eta
u^{m,q}_{N+1}$ are supposed to be known, the
function $\phi^q_N$ is entirely determined.
Observe that if $q=1$, $\partial_\eta u^{m,1}_{N+2}|_{\eta=0}$ is unknown since $\partial_n u^{c,1}_{N}$ is not yet determined, while
as soon as $q\geq 2$, $\partial_\eta u^{m,q}_{N+2}|_{\eta=0}$ is known.

Using transmission condition~\eqref{CTDlcytins} , we infer the following equality satisfied by $u^{m,1}_{N}$ in $\eta=0$:
\begin{align*}-\partial^2_\theta u^{m,1}_N|_{\eta=0} +\beta\partial_n u^{c,1}_{N}\compo\Phi_0=&\phi^1_N
-\int^1_0(\eta-1)\partial^2_\theta\partial_\eta u^{m,1}_N\xdif \eta,\end{align*}
hence the boundary condition imposed to $u^{c,1}_{N}$ on $\partial\mathcal{O}$:
\begin{align*}
\beta\left.\partial_n u^{c,1}_{N}\right|_{\partial\mathcal{O}}-\left.\partial^2_\theta u^{c,1}_{N}\right|_{\partial\mathcal{O}}
=&\Biggl(\phi^1_N-\int^1_0(\eta-1)\partial^2_\theta\partial_\eta u^{m,1}_N\xdif \eta\Biggr)\compo\Phi^{-1}_0.
\end{align*}
Thus the function $u^{c,1}_{N}$ is solution to the following problem:
\begin{equation}
\left\{\begin{aligned}
&\Delta u^{c,1}_{N}=0,\text{ in $\mathcal{O}$,}\\
&-\left.\partial^2_\theta u^{c,1}_{N}\right|_{\partial\mathcal{O}} +\beta\left.\partial_n u^{c,1}_{N}\right|_{\partial\mathcal{O}}
=\Biggl(\phi^1_N-\int^1_0(\eta-1)\partial^2_\theta\partial_\eta u^{m,1}_N\xdif \eta\Biggr)\compo\Phi^{-1}_0,\\
&\int_{\partial\mathcal{O}}u^{c,1}_{N}\xdif_{\partial\mathcal{O}}=0.
\end{aligned}\right.\label{Vc-1}
\end{equation}
{In the membrane $u^{m,1}_N$ is defined by }
\begin{equation}
u^{m,1}_N=\int^s_0\partial_\eta u^{m,q}_{N}\xdif\eta+u^{c,q}_N\compo\Phi_0.\label{Vm-1}
\end{equation}
If $q\geq2$, $u^{m,q}_N|_{\eta=1}$ is entirely determined by the equality:
\begin{align*}-\partial^2_\theta u^{m,q}_N|_{\eta=1} 
=&\beta\partial_n u^{c,q}_{N+1-q}\compo\Phi_0+\phi^q_N-\int^1_0(\eta-1)\partial^2_\theta\partial_\eta u^{m,q}_N\xdif \eta,\intertext{hence}
u^{m,q}_N(s,\theta)&=\int^s_1\partial_\eta u^{m,q}_{N}\xdif\eta+u^{m,q}_N|_{\eta=1}.
\end{align*}
The potential $u^{c,q}_N$ satisfies the following boundary value problem:
\begin{equation}
\left\{\begin{aligned}
&\Delta u^{c,q}_{N}=0,\text{ in $\mathcal{O}$,}\\
&\left. u^{c,q}_{N}\right|_{\partial\mathcal{O}}=u^{m,q}_{N}\compo\Phi^{-1}_0 .
\end{aligned}\right.\label{uc20}
\end{equation}
Observe that for $q\geq1$, $\partial_\eta u^{m,q}_{N+2}$ is then entirely determined by:
\begin{align*}
\partial_\eta u^{m,q}_{N+2}&=\int^s_1\Biggl(-\kappa\left(3\eta\partial^2_\eta u^{m,q}_{N+1}
+\partial_\eta u^{m,q}_{N+1}\right)\\&-\partial^2_\theta u^{m,q}_{N}-\eta\kappa\partial^2_\theta u^{m,q}_{N-1}+\eta\kappa'
\partial_\theta u^{m,q}_{N-1}\Biggr)\xdif\eta.\end{align*}
Therefore, we have proved that for all $N\geq-1$, for $q\in\N^*$, the
functions $u^{c,q}_N$ and $u^{m,q}_N$ are uniquely determined. 
\begin{rmrk}[Regularity]
 Observe
that these functions are the potentials given in Theorem~\ref{thm2}.
We leave the reader verify by induction that the following regularities hold.
Let $q\in\N^*$, $N\geq-1$ and $p\geq 1$.
Let $\phi$ belong to $\sob[N+p-3/2]{\partial\W_h}$.
\begin{subequations}
 \begin{align}
u^{c,q}_{-1}&\in\sob[1+N+p]{\mathcal{O}},\notag\\
u^{m,q}_{-1}&\in\Cscr^\infty\left([0,1];\sob[1/2+N+p]{\Tr}\right),\notag\\
\forall k=0,\cdots, N,\notag\\
u^{c,q}_{k}&\in\sob[1+N+p-2{[k/2]}]{\mathcal{O}},\\
u^{m,q}_k&\in\Cscr^\infty\left([0,1];\sob[1/2+N+p-2{[(k+1)/2]}]{\Tr}\right).
\end{align}\label{regulariteD}
\end{subequations}
Moreover, there exists a constant $C_{N,\mathcal{O},p}$ independant on $h$ and $\beta$ such that:
\begin{subequations}
\begin{align}
\sup_{\eta\in[0,1]}\left\|u^{m,q}_{-1}(\eta,\cdot)\right\|_{\sob[1/2+N+p]{\Tr}}&\leq C_{N,\mathcal{O},p}\|\frm\|_{\sob[N+p-3/2]{\partial\mathcal{O}}},\\
\left\|u^c_{-1}\right\|_{\sob[1+N+p]{\mathcal{O}}}&\leq C_{N,\mathcal{O},p}\|\frm\|_{\sob[N+p-3/2]{\partial\mathcal{O}}},\\
\forall k=0,\cdots, N,\notag\\
\sup_{\eta\in[0,1]}\left\|u^{m,q}_k(\eta,\cdot)\right\|_{\sob[1/2+N+p-2{[(k+1)/2]}]{\Tr}}&\leq C_{N,\mathcal{O},p}\|\frm\|_{\sob[N+p-3/2]{\partial\mathcal{O}}},\\
\left\|u^{c,q}_k\right\|_{\sob[1+N+p-{[k/2]}]{\mathcal{O}}}&\leq C_{N,\mathcal{O},p}\|\frm\|_{\sob[N+p-3/2]{\partial\mathcal{O}}}.
\end{align}\label{estimVmkVckD}
\end{subequations}
\end{rmrk}
\subsection{Error estimates of Theorem~\ref{thm2}}
Let us now prove Theorem~\ref{thm2}.
Let $q\in\N^*$ and $N\in\N$.
The complex parameter  $\alpha$ satisfies \eqref{alpha} with Hypothesis~\ref{hypo1}.
Let $\phi$ belong to $\sob[N+3/2+q]{\partial \W_h}$. 
Let $r^{c,q}_N$ and $r^{m,q}_N$ be the functions defined by:
 \begin{align*}
\begin{cases} r^{c,q}_N=u-\sum^N_{k=-1}u^{c,q}_kh^k,\,\text{in $\mathcal{O}$},\\
 r^{m,q}_N=u\compo\Phi-\sum^N_{k=-1}u^{m,q}_kh^k,\,\text{in $\Cyl$}.
 \end{cases}
 \end{align*}
We have to prove that there exists a constant $C_{\mathcal{O},N}>0$ depending only on the domain $\mathcal{O}$ and on $N$ such that 
\begin{subequations}
 \begin{align}
\|r^{c,q}_N\|_{\sob[1]{\mathcal{O}}}&\leq C_{\mathcal{O},N}\|\frm\|_{\sob[N+3/2+q]{\partial\mathcal{O}}} \max\left(\sqrt{\frac{h}{|\beta|}},\sqrt{h}\right)h^{N+1/2},\\
\|r^{m,q}_N\|_{H^1_\gfrak\left(\Cyl\right)}&\leq C_{\mathcal{O},N}\|\frm\|_{\sob[N+3/2]{\partial\mathcal{O}}}h^{N+1/2}.\label{51b}
 \end{align}\label{hypo22resultat2D}
\end{subequations}
If $\phi$ belongs to $\sob[N+5/2+q]{\partial\W_h}$, we have
$$\left\| r^{c,q}_N\right\|_{\sob[1]{\mathcal{O}}}\leq
 C_{\mathcal{O},N}\|f\|_{\sob[N+5/2+q]{\Tr}}h^{N+1}.$$
\begin{proof}
The proof of Theorem~\ref{thm2} is similar to the proof of Theorem~\ref{thm1}. 
Since $\phi$ belongs to $\sob[N-1/2]{\partial \W_h}$, according to the
previous lemma, the couples of  functions 
$\left(r^{c,q}_{N},r^{m,q}_{N}\right)$ and  $\left(r^{c,q}_{N+1},r^{m,q}_{N+1}\right)$
 are well defined and belong to 
$\sob[1]{\mathcal{O}}\times H^1_\gfrak\left(\Cyl\right)$.

Denote by $\tilde{g}_N$ the following function defined on $\Cyl$:
\begin{align}
\begin{split}\tilde{g}_N=&\kappa\left(3\eta\partial^2_\eta u^{m,q}_{N}+\partial_\eta u^{m,q}_{N}\right)
+3\eta^2\kappa^2\partial^2_\eta u^{m,q}_{N-1}+2\eta\kappa^2\partial_\eta u^{m,q}_{N-1}+\partial^2_\theta u^{m,q}_{N-1}\\&
+\eta^3\kappa^3\partial^2_\eta u^{m,q}_{N-2}+\eta^2\kappa^3\partial_\eta u^{m,q}_{N-2}+\eta\kappa\partial^2_\theta u^{m,q}_{N-2}
-\eta\kappa'\partial_\theta u^{m,q}_{N-2}\\&+
h\Biggl(3\eta^2\kappa^2\partial^2_\eta u^{m,q}_{N}+2\eta\kappa^2\partial_\eta u^{m,q}_{N}+\partial^2_\theta u^{m,q}_{N}
\\&+\eta^3\kappa^3\partial^2_\eta u^{m,q}_{N-1}+\eta^2\kappa^3\partial_\eta u^{m,q}_{N-1}+\eta\kappa\partial^2_\theta u^{m,q}_{N-1}
-\eta\kappa'\partial_\theta u^{m,q}_{N-1}
\Biggr)\\&+
h^2\Biggl(\eta^3\kappa^3\partial^2_\eta u^{m,q}_{N}+\eta^2\kappa^3\partial_\eta u^{m,q}_{N}+\eta\kappa\partial^2_\theta u^{m,q}_{N}
-\eta\kappa'\partial_\theta u^{m,q}_{N}
\Biggr)\end{split}
\end{align}
According to the previous lemma and since $\phi$ belongs to $\sob[N-1/2]{\partial\W_h}$, the above function $\tilde{g}_N$ belongs to 
$\Cscr^\infty\left([0,1];\sob[-1/2]{\Tr}\right)$ and the function $\partial_\eta V^{m,q}_N$ belongs to 
$\Cscr^\infty\left([0,1];\sob[3/2]{\Tr}\right)$.  Moreover, there exists a constant $C_{N,\mathcal{O}}$ such that
\begin{align}
\begin{cases}
\sup_{\eta\in[0,1]}\|\tilde{g}_N(\eta,\cdot)\|_{H^{-1/2}\left(\Tr\right)}\leq C_{N,\mathcal{O}}\|f\|_{\sob[N-1/2]{\Tr}},\\
\sup_{\eta\in[0,1]}\left\|\partial_\eta u^{m,q}_N(\eta,\cdot)\right\|_{H^{3/2}\left(\Tr\right)}\leq C_{N,\mathcal{O}}\|f\|_{\sob[N-1/2]{\Tr}}.\end{cases}
\label{hypo22estimgNcytins}
\end{align}
The functions $r^{c,q}_N$ and $r^{m,q}_N$ satisfy the following problem:
\begin{align*}
&\Delta r^{c,q}_N=0,\text{ in $\mathcal{O}$},\\
&\partial_\eta\left(\frac{1+h\eta\kappa}{h}\partial_\eta r^{m,q}_N\right)+\partial_\theta\left(\frac{h}{1+h\eta\kappa}\partial_\theta r^{m,q}_N\right)=\frac{-h^N}{(1+h\eta\kappa)}\tilde{g}_N,
\intertext{with transmission conditions:}
&\beta h^{1+q}\partial_nr^{c,q}_N\compo\Phi_0=\frac{1}{h}\left(\left.\partial_\eta r^{m,q}_N\right|_{\eta=0}+\beta h^{N+1+q}\left(\partial_n u^{c,q}_{N-1}\compo\Phi_0
+h\partial_\eta u^{m,q}_{N}\compo\Phi_0\right)\right),\\
&r^{c,q}_N\compo\Phi_0=\left.r^{m,q}_N\right|_{\eta=0},
\intertext{with boundary condition}
&\left.\partial_\eta r^{m,q}_N\right|_{\eta=1}=0,
\intertext{and with gauge condition}
&\int_{\partial\mathcal{O}}r^{c,q}_N\,\xdif\sigma=0.
\end{align*}
By multiplying the above equality by $\overline{r_N}$ and by integration by parts, we infer that:
\begin{align}\begin{split}\beta h^{1+q}\left\|\xdif r^{c,q}_N\right\|^2_{\Lambda^1\leb[2]{\mathcal{O}}}&+\left\|\xdif r^{m,q}_N\right\|^2_{\Lambda^1L^2_\gfrak\left(\Cyl\right)}
=- h^N\int_{\Cyl}\tilde{g}_N(\eta,\theta)\overline{r^{m,q}_N}(\eta,\theta)\,\xdif\eta\,\xdif\theta\\
&+\beta h^{N+1+q}\int_{\Tr}\left(\partial_\eta u^{c,q}_{N-1}\compo\Phi_0+h\partial_\eta u^{c,q}_N\compo\Phi_0\right)\left.\overline{r^{m,q}_N}\right|_{\eta=0}\,\xdif\theta.\end{split}\label{hypo22equalitycytins}
 \end{align} 
The end of the proof is similar to Theorem~\ref{thm1}.
Using the positivity of $\Re(\beta)$ we straight infer estimate \eqref{51b}
of $r^{m,q}_N$.
To obtain the estimates of 
$r^{c,q}_N$, we write:
$$r^{c,q}_N=r^{c,q}_{N+q}+\sum^q_{k=1} u^{c,q}_{N+k}h^{N+k}.$$
\end{proof}

\subsection{The case $\alpha=o(h^N)$, $\forall N\in\N$}
Now, we suppose that Hypothesis~\ref{hypo2} holds. In this case, we prove that $u^{c}$ and $u^{m}$ may be approximated by $U^{c}$ and $U^{m}$, which are solution to:
\begin{subequations}
\begin{align}
&\Delta U^{m}=0,\text{ in $\mathcal{O}_h$},\\
&\left.\partial_\eta U^{m}\right|_{\partial\mathcal{O}}=0,\quad\left.\partial_\eta U^{m}\right|_{\partial\W_h}=\phi,\\
&\int_{\partial\mathcal{O}}U^{m} \xdif\sigma=0.
\end{align}\label{hypo3Um}
\end{subequations}
and 
\begin{subequations}
\begin{align}
&\Delta U^{c}=0,\text{ in $\mathcal{O}$},\\
&\left.U^{c}\right|_{\partial\mathcal{O}}=\left.U^{m}\right|_{\partial\mathcal{O}}.
\end{align}\label{hypo3Uc}
\end{subequations}
Actually, we have the following lemma:
\begin{lemm}
Let $\phi$ belong to $\sob[-1/2]{\partial\W_h}$. Let $(u^{c},u^{m})$ be the solution to Problem~\eqref{diel11}, and $U^{m}$ and $U^{c}$ be defined respectively by \eqref{hypo3Um} and \eqref{hypo3Uc}.
Then, we have:
\begin{align}
&\left\|u^{m}-U^{m}\right\|_{\sob[1]{\mathcal{O}_h}}\leq C_{\mathcal{O}}|\alpha|\left|\phi\right|_{\sob[-1/2]{\partial\W_h}},\\
&\left\|u^{c}-U^{c}\right\|_{\sob[1]{\mathcal{O}}}\leq C_{\mathcal{O}}\sqrt{|\alpha|}\left|\phi\right|_{\sob[-1/2]{\partial\W_h}}.
\end{align}
\end{lemm}
\begin{proof}
Denote by $w^{c}$ and $w^{m}$ the following functions:
$$w^{c}=u^{c}-U^{c},\quad w^{m}=u^{m}-U^{m},$$
and let $\phi$ belong to $\sob[-1/2]{\partial\W_h}$. We have:
\begin{subequations}
\begin{align}
&\Delta w^{c}=0,\text{ in $\mathcal{O}$},\\
&\Delta w^{m}=0,\text{ in $\mathcal{O}_h$},\\
&\alpha\left.\partial_nw^{c}\right|_{\partial\mathcal{O}}=\left.\partial_nw^{m}\right|_{\partial\mathcal{O}}-\alpha \left.\partial_nU^{c}\right|_{\partial\mathcal{O}},\\
&\left.w^{c}\right|_{\partial\mathcal{O}}=\left.w^{m}\right|_{\partial\mathcal{O}},\\
&\left.\partial_\eta w^{m}\right|_{\partial\W_h}=0,\\
&\int_{\partial\mathcal{O}}w^{m} \xdif\sigma=0.
\end{align}
\end{subequations}
Thus we infer:
\begin{align}
\alpha\int_{\mathcal{O}}\left|\nabla w^{c}\right|^2\dvol_{\mathcal{O}}+\int_{\mathcal{O}_h}\left|\nabla w^{m}\right|^2\dvol_{\mathcal{O}_h}=
\alpha\int_{\partial\mathcal{O}}\left.\partial_nU^{c}\right|_{\partial\mathcal{O}}\overline{w^{m}}\xdif\sigma.
\end{align}
It is well-known that :
$$\left\|U^{m}\right\|_{\sob[1]{\mathcal{O}_h}}\leq C_{\mathcal{O}}\left|\phi\right|_{\sob[-1/2]{\partial\W_h}},$$
and 
$$\left\|U^{c}\right\|_{\sob[1]{\mathcal{O}_h}}\leq C_{\mathcal{O}}\left|U^{m}|_{\partial\mathcal{O}}\right|_{\sob[1/2]{\partial\mathcal{O}}}.$$
Since $\alpha$ satisfies \eqref{alpha} we infer,  $$\left\|w^{m}\right\|_{\sob[1]{\mathcal{O}_h}}\leq C_{\mathcal{O}}|\alpha|\left|\phi\right|_{\sob[-1/2]{\partial\W_h}},$$
and thereby $$\left\|w^{c}\right\|_{\sob[1]{\mathcal{O}}}\leq C_{\mathcal{O}}\sqrt{|\alpha|}\left|\phi\right|_{\sob[-1/2]{\partial\W_h}}.$$
\end{proof}
It remains to derive  asymptotics of $U^{m}$ and then these of $U^{c}$. They are similar to asymptotics of $u^{m,q}$ for $q\geq2$: 
we just have to replace $\beta$ by zero. We think the reader may
easily derive these asymptotics from our previous results.

\section*{Conclusion}
In this paper, we have studied the steady state
potentials in a highly contrasted domain with thin layer when Neumann
boundary condition is imposed on the exterior boundary. 
We derived rigorous asymptotics with respect to the thickness of the
potentials in each domain and we gave error estimate in terms of
appropriate Sobolev norm of the boundary data, electromagnetic
parameters of our domain and a constant depending only on the geometry
of the domain.  
It has to be mentionned that for an insulating inner domain (or
equivalently a conducting membrane), the asymptotic expansions start
at the order -1 and mixed or Dirichlet boundary conditions has to be
imposed on the asymptotic terms of the inner domain.

To illustrate these asymptotics, numerical simulations using FEM are forthcomig work
with Patrick Dular from Universit\'e de Li\`ege and Ronan Perrussel from Amp\`ere laboratory of Lyon. 
Few results have been shown at the conference NUMELEC \cite{numelec06} with GetDP\cite{getdp}. The main difficulty in 
illustrating the convergence of our asymptotic consists in the geometrical approximation of the domain: high-order geometric 
elements seem to be necessary.

 \section*{Appendix}

 Let $\star$ denote the Hodge star operator, which maps $0$-forms to $2$-forms, $1$-forms to $1$-forms and $2$-forms to $0$-forms (see Flanders \cite{flanders}).
 We give explicit formulae for the operators $\xdif$,  $\delta$, $\ext$ and $\inte$. 
 These formulae are straightforward consequences of the definition of the operators ${\star}$, $\xdif$ and $\delta={\star}^{-1}\xdif\,{\star}$. 
 We refer the reader to Dubrovin,  Fomenko and Novikov~\cite{dubrovin}. 

 We consider the metric given by the following matrix  $G$ 
 \begin{equation}
 G=\left(\begin{matrix}
   g_{11} & g_{12} \\
   g_{12} & g_{22} \\
 \end{matrix}\right).
 \end{equation}
 We denote by $|G|$ the determinant of $G$. The inverse of $G$ is denoted by $G^{-1}$
  $$G^{-1}=(g^{ij})_{ij},$$ and we suppose that the signature of $G$ is equal to 1. Thereby, the operator ${\star}^2$ is equal to $\id$ on the space of 
 0-forms and 2-forms and it is equal to $-\id$ on 1-forms. 
  \subsection{Star operator in  $\Er^2$}
 \subsubsection{On 0-forms and on 2-forms}
 Let $T$  be a 0-form and let $S$ be the  2-form $\nu\,\xdif y^1\xdif y^2$. Then ${\star}T$ is the 2-form $\mu\,\xdif y^1\xdif y^2$
 and ${\star}S$ is the 0-form $f$.
 The following identities hold:
 \begin{align*}
 \mu&=\sqrt{|G|}T,\\
 f&=\frac{1}{\sqrt{|G|}}\nu.
 \end{align*}
 \subsubsection{On 1-forms}
 Let  $T$ be the 1-form $T_1\,\xdif y^1+T_2\,\xdif y^2$. Then ${\star}T$
 is the 1-form $\mu_1\,\xdif y^1+\mu_2\,\xdif y^2$, and we have the following formulae:
 \begin{align*}
 \mu_1&=-\sqrt{|G|}\left(g^{12}T_1+g^{22}T_2\right),\\
 \mu_2&=\sqrt{|G|}\left(g^{11}T_1+g^{12}T_2\right).
 \end{align*}

 \subsection{The action of $\xdif$ acting on 0-forms in $\Er^2$}
 Let  $\mu$ be a 0 form, then $\xdif \mu$ has the following expression:
 \begin{align*}
 &\xdif \mu=\frac{\partial \mu}{\partial y^1}\xdif y^1+\frac{\partial
   \mu}{\partial y^2}\xdif y^2.
 \end{align*}

 \subsection{The action of $\xdelta$ acting on 1-forms on $\Er^2$}
 Let $\mu$ be the 1-form $\mu_{1}\xdif y^1+\mu_2\xdif y^2$, and define $\xdelta\mu=\alpha$.
 The 0-form $\alpha$ is equal to:
 \begin{align*}
 \alpha=-\frac{1}{\sqrt {|G|}}
 \biggl\{&\frac{\partial}{\partial y_1}\Bigl(\sqrt{|G|}\left(
 g^{11}\mu_1+g^{12}\mu_2\right)\Bigr)\\
 &+\frac{\partial}{\partial y_2}\Bigl(\sqrt{|G|}\left(
 g^{12}\mu_1+g^{22}\mu_2\right)\Bigr)\biggr\}.
 \end{align*}
 \subsection{The exterior product of a 1-form with a 0-form}
 Let $N$ be the 1-form $N_1\xdif y^1+N_2\xdif y^2$ and $f$ be a 0-form. The exterior product of $\ext(N)f$ is:  
 \begin{align*}
 \ext(N)f=f N_1\xdif y^1+f N_2\xdif y^2.
 \end{align*}

 \subsection{The interior product of a 1-form with a 1-form}
 Let $N$ and $\mu$ be the 1-forms $N_1\xdif y^1+N_2\xdif y^2$, and 
 $\mu_{1}\xdif y^1+\mu_2\xdif y^2$.
 Then 0-form $\inte(N)\mu$ has the following expression:
 \begin{align*}
 \inte(N)\mu=N_1\left(\mu_1g^{11}+\mu_2g^{12}\right)+
 N_2\left(\mu_1g^{12}+\mu_2g^{22}\right).
 \end{align*}

\bibliographystyle{plain}
\bibliography{bibli_asympt}
\end{document}